\documentclass[final,onefignum,onetabnum]{siamart171218}

\usepackage{lipsum}
\usepackage{amsfonts}
\usepackage{amsmath}
\usepackage{amssymb}
\usepackage{dsfont}
\usepackage{graphicx}
\usepackage[caption=false]{subfig}
\usepackage{epstopdf}
\usepackage{algorithmic}
\usepackage{amsopn}
\usepackage{siunitx}

\newcommand{\R}{\mathbb{R}}
\newcommand{\J}{\mathcal{J}}
\newcommand{\E}{\mathcal{E}}
\newcommand{\norm}[1]{\Vert #1 \Vert}
\newcommand{\sqnorm}[1]{\Vert #1 \Vert_2^2}

\newcommand{\argmin}{\text{argmin}}
\newcommand{\argmaxsub}[1]{\underset{{ #1 }}{{\rm argmax}}}
\newcommand{\argminsub}[1]{\underset{{ #1 }}{{\rm argmin}}}

\newcommand{\prox}{{\rm prox}}
\newcommand{\diag}{{\rm diag}}
\newcommand{\NabFD}{\nabla^{+}}

\newcommand{\NabCD}{\nabla^{\pm}}
\newcommand{\parFD}{\partial^{+}}
\newcommand{\parCD}{\partial^{\pm}}
\newcommand{\laplaceFD}{\Delta^{+}}

\newcommand{\termabb}[2]{\emph{#1} (\emph{#2})}
\newcommand{\dontshow}[1]{}

\DeclareMathOperator*{\mydef}{\mathrel{\mathop:}=}

\ifpdf
  \DeclareGraphicsExtensions{.eps,.pdf,.png,.jpg}
\else
  \DeclareGraphicsExtensions{.eps}
\fi

\newsiamremark{remark}{Remark}
\newsiamremark{hypothesis}{Hypothesis}
\crefname{hypothesis}{Hypothesis}{Hypotheses}
\newsiamthm{claim}{Claim}

\headers{Enhancing Compressed 4D PAT by Simultaneous Motion Estimation}{F. Lucka, N. Huynh, M. Betcke, E. Zhang, P. Beard, B. Cox, S. Arridge}

\title{Enhancing Compressed Sensing 4D Photoacoustic Tomography by Simultaneous Motion Estimation\thanks{
\funding{This work was supported in parts by the Engineering and Physical Sciences Research Council, UK (EP/K009745/1), the European Union project FAMOS (FP7 ICT, Contract 317744), the European Union's Horizon 2020 research and innovation programme H2020 ICT 2016-2017 under grant agreement No 732411 (as an initiative of the Photonics Public Private Partnership), the Netherlands Organisation for Scientific Research (NWO 613.009.106/2383) and the National Institute of General Medical Sciences of the National Institutes of Health under grant number P41 GM103545-18.}}}

\author{Felix Lucka\thanks{Computational Imaging, Centrum Wiskunde \& Informatica (CWI), Science Park 123, 1098 XG Amsterdam, The Netherlands and Department of Computer Science, University College London, WC1E 6BT London, UK (\email{Felix.Lucka@cwi.nl}, \url{felixlucka.github.io/}).}
\and Nam Huynh\thanks{Department of Medical Physics and Bioengineering, University College London, WC1E 6BT London, UK}
\and Marta Betcke\thanks{Department of Computer Science, University College London, WC1E 6BT London, UK.}
\and Edward Zhang\footnotemark[3]
\and Paul Beard\footnotemark[3]
\and Ben Cox\footnotemark[3]
\and Simon Arridge\footnotemark[4]
}

\ifpdf
\hypersetup{
  pdftitle={Enhancing Compressed Sensing 4D PAT by Simultaneous Motion Estimation},
  pdfauthor={F. Lucka}
}
\fi

\begin{document}

\maketitle

\begin{abstract}
A crucial limitation of current high-resolution 3D photoacoustic tomography (PAT) devices that employ sequential scanning is their long acquisition time. In previous work, we demonstrated how to use compressed sensing techniques to improve upon this: images with good spatial resolution and contrast can be obtained from suitably sub-sampled PAT data acquired by novel acoustic scanning systems if sparsity-constrained image reconstruction techniques such as total variation regularization are used. Now, we show how a further increase of image quality can be achieved for imaging dynamic processes in living tissue (4D PAT). The key idea is to exploit the additional temporal redundancy of the data by coupling the previously used spatial image reconstruction models with sparsity-constrained motion estimation models. While simulated data from a two-dimensional numerical phantom will be used to illustrate the main properties of this recently developed joint-image-reconstruction-and-motion-estimation framework, measured data from a dynamic experimental phantom will also be used to demonstrate its potential for challenging, large-scale, real-world, three-dimensional scenarios. The latter only becomes feasible if a carefully designed combination of tailored optimization schemes is employed, which we describe and examine in more detail. 
\end{abstract}

\begin{keywords}
Photoacoustic tomography, dynamic imaging, compressed sensing, simultaneous motion estimation, variational regularization.
\end{keywords}

\begin{AMS}
92C55, 65R32, 94A08, 94A12, 65K10.
\end{AMS}

\section{Introduction} \label{sec:Intro}

\subsection{Compressed Sensing Photoacoustic Tomography} \label{subsec:IntroPAT}

Optical absorption of biological tissues is a desirable source of image contrast for a variety of clinical and preclinical applications. In particular, its wavelength dependence provides spectroscopic (chemical) information on the absorbing molecules (chromophores). \termabb{Photoacoustic Tomography}{PAT} is an "Imaging from Coupled Physics"-technique \cite{ArSc12} that employs laser-generated ultrasound (US) to obtain optical absorption images with the high spatial resolution of US. For recent reviews on the physical principles, technical realizations and (pre-)clinical applications of PAT, we refer the reader to \cite{Wa09,Bea11,NiCh14,ZhYaWa2016}. \\
In \cite{ArBeBeCoHuLuOgZh16}, we discussed the particular challenges of acquiring high quality three-dimensional (3D) photoacoustic (PA) images with sequential scanning schemes, such as the \termabb{Fabry-P\'{e}rot based PA scanner}{FP scanner}:
To reach a spatial resolution less than one hundred \si{\micro \meter}, acoustic waves containing frequencies up to a few tens of \si{\mega \hertz} have to be sampled over \si{\centi \meter} scale apertures. For a scanning pattern to satisfy the spatial Nyquist criterion, sampling intervals in the order of tens of \si{\micro \meter} have to be chosen, which leads to several thousand detection points and thereby, long acquisition times. This imposes a severe limit for dynamic PAT (4D PAT), i.e., imaging dynamic anatomical and physiological events in high resolution in real time, an area of research of increasing interest \cite{DeGoMcShRa17}. The key observation to overcome this limitation is that the Nyquist criterion is often too conservative because it guarantees perfect recovery of the broad class of images that are band-limited but otherwise arbitray. However, images of absorbing tissue structures come from a much smaller sub-class of images, as they typically also have a rather low spatial complexity (or a high \textit{sparsity}). Therefore, data recorded in a conventional, regularly sampled fashion, satisfying the Nyquist criterion, is often highly redundant. \termabb{Compressed Sensing}{CS} \cite{CaRoTa06,Do06,FoRa13} techniques exploit this fact by combining sub-sampling schemes that try to maximize the non-redundancy of the data with image reconstruction approaches that employ sparsity-constraints. In \cite{ArBeBeCoHuLuOgZh16}, we demonstrated the implementation of CS techniques to accelerate 3D PAT acquisition by using spatial sparsity constraints. In the context of 4D PAT, such techniques can be employed to reconstruct each temporal frame separately, i.e., as a \termabb{frame-by-frame}{fbf} image reconstruction method.

\subsection{Spatio-Temporal Image Reconstruction} \label{subsec:IntroSTImaging}

In this work, we show that another significant acceleration can be obtained by also accounting for the temporal evolution of the target within a full spatio-temporal reconstruction scheme. A wide range of such approaches have been proposed for different applications and dynamics. If the dynamics between separate frames are sufficiently simple (e.g., affine deformations), low-dimensional parametric models can often be used to efficiently constrain the image reconstruction in time. An application to PAT is demonstrated in \cite{ChNg16}  and theoretical analysis of such approaches can be found in \cite{HaQu16,Ha16,Ha15}. In such situations, the aim is often rather to compensate for the motion (see, e.g., \cite{McHaScKi13} for an overview on compensating respiratory motion) than to resolve it, which is our main aim here. 
Several approaches rely on extending popular spatial constraints into time. Incorporating $\ell_2$ regularization of the temporal differences between frames is examined in \cite{ScLo02,ScLoWoVa02} and recently, extending $\ell_1$ functionals such as total variation functional and its higher order variants to spatio-temporal settings have been proposed and have been shown to work well for certain dynamics, e.g.,  \cite{HoKu14,ScHoScBrSt17}. 
In the Bayesian approach to inverse imaging problems, spatio-temporal methods are commonly refered to as \textit{Kalman filtering or smoothing}: Filtering refers to reconstructing each image frame based only on measured data up to that point in time, most often done via updating the previous image frame based on the most recent data. While this is the only option for real-time or \textit{online} image reconstruction, it is also popular in \textit{offline} image reconstruction due to its lower computational complexity compared to smoothing, which refers to estimating each image frame based on the whole set of measured data. See Section 4 in \cite{KaSo05} for a general introduction and further references to Kalman fitering and \cite{SoCuHaMa16} for recent work on this topic. 
In the context of compressed sensing applications, low-rank-type models have been examined extensively, see, e.g., \cite{HaLi10,WaXiLiWaAn14,TrDiAtAr14,RaMoNaFe17}. These models rely on strong spatio-temporal decomposition assumptions which are very effective when fulfilled but not appropriate for every dynamics. \\
In this work, we adopt a very general spatio-temporal modelling framework introduced in \cite{BuDiSc16} that can encode a-priori information about a wide range of dynamics: It formulates an explicit PDE model for the image dynamics and then jointly estimates the image sequence and the corresponding motion field by minimizing a variational energy. An overview of similar approaches to joint image reconstruction and motion estimation can be found in the introduction of \cite{BuDiSc16}, which also contains theoretical analysis of this approach. While it was used for 2D dynamic computed tomography reconstruction in \cite{BuDiFrHaHeSi17}, we present the first application to a challenging, large-scale 3D dynamic problem with experimental data, which  also requires the development of tailored numerical optimization schemes. 

\subsection{Structure} \label{subsec:Structure}

The remainder of the paper is organized as follows: Section \ref{sec:Background} introduces the mathematical modeling of dynamic PAT and illustrated the limitations of reconstruction approaches that only account for spatial sparsity. Based on this, a variational spatio-temporal image reconstruction framework based on joint motion estimation is presented in Section \ref{sec:JointImRecMotionEst}. Section \ref{sec:Opt} discusses the numerical solution of the optimization problems that originate from the variational approach and in Section \ref{sec:Res}, we present results with a simple 2D scenario with simulated data and a challenging 3D scenario with experimental data. Finally, we discuss the results of our work and point to future directions of research in Section \ref{sec:DisOutCon}. Table \ref{tbl:Abb} lists all commonly occurring abbreviations for reference.

\begin{table}
{\footnotesize
\caption{\label{tbl:Abb} List of commonly occurring abbreviations.} 
\begin{center}
\begin{tabular}{@{}lll}        
\bf Abbreviation & \bf Meaning & \bf Reference  \\ \hline
ACS     & alternate convex search & Sec. \ref{sec:BiConv} \\
ADMM & alternating direction method of multipliers & Sec. \ref{sec:ADMM}, Alg. \ref{algo:ADMM} \\
fbf & frame-by-frame & Sec. \ref{subsec:PreWork}  \\
FP  & Fabry-P\'{e}rot  & Sec. \ref{subsec:IntroPAT} \\
mIP & maximum intensity projection & Fig. \ref{fig:LimFBF} \\
NNLS & non-negative least squares & Sec. \ref{sec:NumPhan} \\
(Q)PAT & (quantitative) photoacoustic tomography & Sec. \ref{sec:Intro} \\
PDHG & primal dual hybrid gradient & Sec \ref{sec:PDHG}, Alg. \ref{algo:PDHG} \\
TV  & total variation regularization & Sec. \ref{sec:OptFlow} \\
TVTVL2  & Joint image reconstruction and motion estimation approach  & Sec. \ref{sec:OptFlow},  \eqref{eq:TVTVL2SemiCont} \\ \hline
\end{tabular}
\end{center}
}
\end{table}

\section{Background and Previous Work} \label{sec:Background}

\subsection{Sequential Acquisition  of Compressed Dynamic PAT} \label{subsec:DynPAT}

\begin{figure}[t]
   \centering
\includegraphics[width=\textwidth]{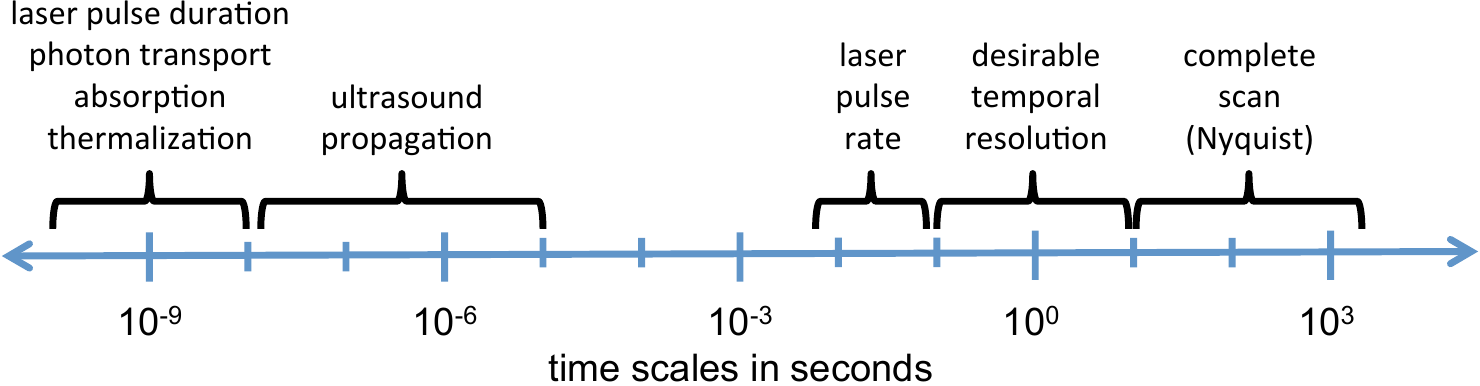}
\caption{Sketch of the relevant time scales in high resolution 4D PAT with sequential acquisition: PAT is particularly suited to image dynamic processes in living tissue that are related to blood circulation. For this, one would ideally like to obtain a temporal resolution close to the heart beat ($\sim1$\si{\second} for humans, $\sim0.1$\si{\second} for mice). As the processes that contribute to a single PA signal take place within nanoseconds to microseconds, the main temporal limitation of sequential acquisition systems is given by the excitation laser pulse rate. Lasers with sufficiently high pulse energies are currently limited to $\sim200$\si{\hertz}, which typically leads to longer acquisition times for a complete scan than what is desired. For instance, scanning $20\,000$ locations with a $20$\si{\hertz} laser takes $1000$\si{\second} while using a  $200$\si{\hertz} laser in combination with a multi-beam read-out system as described in \cite{HuOgZhCoBe16} takes $12.5$\si{\second}. In the latter case, applying compressed sensing with a sufficiently high sub-sampling factor would yield the desired temporal resolution.}
   \label{fig:timescalesPAT}
\end{figure}

Let us denote the biological tissue to be imaged by $\Omega \in \R^d$ ($d=2,3$), the space variable by $r \in \Omega$, the measurement interval by $[0, \mathcal{T}]$ and the (continous) time variable by $\tau \in [0, \mathcal{T}]$. 
A reasonable mathematical model of dynamic PAT has to make certain assumptions about the different time scales involved in signal generation and measurement, in particular if the PA signal is scanned in a sequential manner. Firstly, as described in more detail in Section 1.1.\ of \cite{ArBeCoLuTr16}, the photoacoustic effect is only significant if the laser pulse duration, photon transport, photon absorption by chromophores and subsequent thermalization take place sufficiently fast, i.e., within a few nanoseconds. The induced, local pressure increase $p: \Omega \rightarrow \R_+$ initiates a broadband acoustic pulse that travels through $\Omega$ within a few microseconds. Therefore, this part of the signal generation is commonly modelled as an initial value problem for the wave equation:
\begin{equation}
(\partial_{\tau \tau} - c^2 \Delta) \breve{p}(r,\tau) = 0 \enspace, \qquad \breve{p}(r,\tau = 0) = p \enspace , \qquad \partial_\tau \breve{p}(r,\tau = 0) = 0 \enspace. \label{eq:PATfwd}
\end{equation}
This approximates the whole optical part as instantaneous, which is equivalent to assuming the tissue remains at rest until the thermalization is complete. Sequential scanning systems can only measure a single spatial projection of $\breve{p}(r,\tau)$ over a sensor surface $\mathcal{S} \subset \partial \Omega$ for each pulse of the excitation laser:
\begin{equation}
f_{m,l} = \int_{[0,\mathcal{T}]} \int_{\mathcal{S}} \breve{p}(r, \tau) \phi_m(r) \psi_{l}(\tau) \; d r \; d \tau, \quad m = 1,\ldots,M, \quad l = 1,\ldots,M_\tau \enspace .
 \label{eq:GenScan}
\end{equation}
where $\phi_m(r)$ describes the spatial window function used for the measurement associated with the $m$-th laser pulse, and $\psi_l$ is the $l$-th temporal window function (we will only consider equidistant temporal point sampling in the following). \\
A single pressure-time series is recorded within a few microseconds, and can therefore be regarded as instantaneous if we are interested in imaging dynamics taking place on the scale of a few seconds or even minutes. However, as described in more detail in \cite{ArBeBeCoHuLuOgZh16}, to form high resolution 3D images, the spatial Nyquist criterion necessitates that several thousand of such time series are recorded. As the pulse repetition rates of conventional excitation lasers are typically limited to tens of \si{\hertz}, this means that the scanning process and the image dynamics interfere - the image is moving while the scanning is taking place - and neglecting this by assuming an instantaneous measurement can lead to severe motion blurring in the reconstructed images. A summary of the relevant time scales is depicted in Figure \ref{fig:timescalesPAT}.  \\
A fully continuous modelling encompassing all the different and interfering spatio-temporal processes described above is of only limited practical  value and will not be pursued here. Instead, we assume that a temporal binning of the sequential acquisitions \eqref{eq:GenScan} into temporal frames, $t=1,\ldots,T$, is chosen in such a way that the initial pressure can be assumed to be static during one frame. We then model the linear mapping of the discretized initial pressure $p_t \in \R^N$ to fully-sampled, discrete data $f_t \in \R^{M M_\tau}$ via \eqref{eq:PATfwd} and \eqref{eq:GenScan} by a time-independent, i.e., instantaneous, operator $A$. In this context, "fully-sampled" refers to an ideal scanning scheme that samples $\mathcal{S}$ as demanded by the spatial Nyquist criterion, although our measurement set-up might practically not allow for doing that within the duration of a single temporal bin. The real measurement is modeled by applying a time-dependent sub-sampling or \emph{compression} operator $C_t \in \R^{M_c M_\tau \times M M_\tau}$ to $f_t$: 
\begin{equation}
f^c_t = C_t f_t = C_t A p_t + \varepsilon_t \enspace, \qquad \qquad t = 1,\ldots,T \enspace, \label{eq:DynPAT}
\end{equation}
where $\varepsilon_t$ accounts for additive measurement noise, which we assume can be modeled as i.i.d.\ standard normal distributed after suitable data pre-processing is carried out. \\
We will mainly use $s$-periodic sequences $C_t$, $t = 1,2,\ldots$, such that for any $t_0 \geqslant 1$,
\begin{equation}
\bar{C}_{t_0} = \begin{pmatrix}
C_{t_0} \\
C_{t_0 + 1} \\
\vdots \\
C_{t_0 + s}
\end{pmatrix} 
\end{equation}
is invertible can be transformed into $\bar{C}_{1}$ by row-permutation. This amounts to splitting a conventional, full scanning pattern $\bar{C} \in \R^{M M_\tau \times M M_\tau}$ consisting of $M$ spatial projections into smaller temporal bins comprising disjoint sub-sets of $M_c$ spatial projections and allows for an intuitive definition of the sub-sampling factor as $M_{{\rm sub}} = M/M_c$. However, the methods presented here can be used for any sequence $\{C_t\}^T_t$.\\
From now on, any reference to time is with respect to the image and measurement dynamics (indexed by $t$), not to the acoustic wave propagation (indexed by $\tau$). Furthermore, we will often ease the notation when dealing with spatio-temporal quantities: Dropping the temporal index $t$ refers to the whole sequence as a vector, e.g., $p \in \R^{N T}$. When spatial operators like the gradient $\nabla$ are applied to such a vectorized dynamic quantity, it is understood as a frame-by-frame application, i.e., $\nabla p$ means $\left( I_T \otimes \nabla \right) p$, where $I_T$ is the $T$ dimensional identity matrix.  

\subsection{Previous Work} \label{subsec:PreWork}

\begin{figure}[htb]
   \centering
\hfill
\subfloat[][$t = 22$, X mIP\label{subfig:LimFBF_22_X}]{\includegraphics[width=0.445\textwidth]{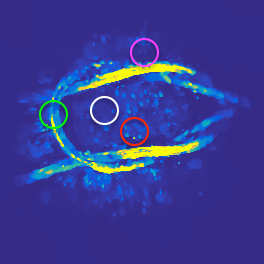}}
\hfill
\subfloat[][$t = 23$, X mIP\label{subfig:LimFBF_23_X}]{\includegraphics[width=0.445\textwidth]{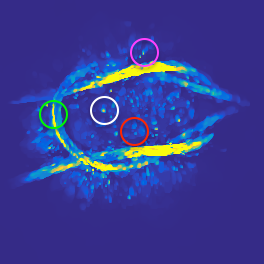}}
\hfill
\\[-10pt]
\hfill
\subfloat[][$t = 22$, Y mIP\label{subfig:LimFBF_22_Y}]{\includegraphics[width=0.445\textwidth]{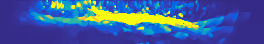}}
\hfill
\subfloat[][$t = 23$, Y mIP\label{subfig:LimFBF_23_Y}]{\includegraphics[width=0.445\textwidth]{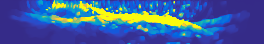}}
\hfill
\\[-10pt]
\hfill
\subfloat[][$t = 22$, Z mIP\label{subfig:LimFBF_22_Z}]{\includegraphics[width=0.445\textwidth]{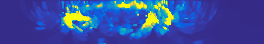}}
\hfill
\subfloat[][$t = 23$, Z mIP\label{subfig:LimFBF_23_Z}]{\includegraphics[width=0.445\textwidth]{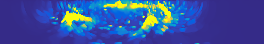}}
\hfill
\\[-10pt]
\hfill
\subfloat[][$t = 22$, X slice\label{subfig:LimFBF_22_sl}]{\includegraphics[width=0.445\textwidth]{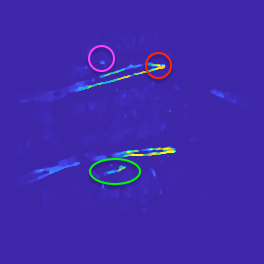}}
\hfill
\subfloat[][$t = 23$, X slice\label{subfig:LimFBF_23_sl}]{\includegraphics[width=0.445\textwidth]{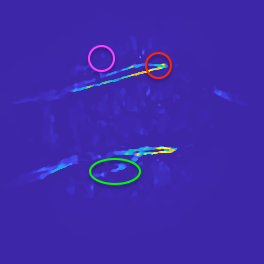}}
\hfill
\\
\caption{Limitations of applying frame-by-frame image reconstruction \eqref{eq:ImRecFBF} to a dynamic PAT data set (details given in Section \ref{sec:ExpPhan}) when using a high sub-sampling factor $M_{{\rm sub}} = 16$: Figures in the left and right column show \termabb{maximum intensity projections}{mIP} along different directions (\protect\subref{subfig:LimFBF_22_X}-\protect\subref{subfig:LimFBF_23_Z})  and slice view visualizations (\protect\subref{subfig:LimFBF_22_sl}-\protect\subref{subfig:LimFBF_23_sl}) of the results $p_{22}$ and $p_{23}$, respectively. One can easily see that image artefacts are not correlated between the two subsequent time frames (same coloured circles in left and right images highlight examples) while the target's motion is.}
   \label{fig:LimFBF}
\end{figure}

In \cite{ArBeBeCoHuLuOgZh16}, we focused on fbf  image reconstruction techniques for \eqref{eq:DynPAT}, i.e., we reconstructed each $p_t$ separately, agnostic to any temporal relationship in the data $f_t^c$. In particular, we showed that variational approaches,
\begin{equation}
\hat{p_t} = \argminsub{p_t \geqslant0} \Bigg\lbrace  \frac{1}{2} \sqnorm{C_t A \, p_t - f_t^c} + \alpha \mathcal{J}(p_t) \Bigg\rbrace \enspace,  \qquad \alpha > 0 \enspace,\label{eq:ImRecFBF}
\end{equation}
that use the regularization functional $\mathcal{J}(p)$ to impose sparsity constraints that encode a-priori knowledge that the images mainly consist of structures of low spatial complexity outperform linear reconstructions such as \emph{time-reversal} or other \emph{back-projection}-type approaches \cite{FiPa04,XuWa04b,ArBeCoLuTr16}. Similar studies by others confirm these results \cite{PrLe09,GuLiSoWa10,ZhWaZh12,MeWaLiSo12,MeWaYiLiSo12,HuWaNiWaAn13,BoLaStGiMaBr17,ElLuZhBeCo17}. As \eqref{eq:ImRecFBF} has to be solved by iterative optimization schemes, fbf image reconstruction is appealing from a computational perspective. However, as it can only encode spatial a-priori information, its ability to obtain good quality images from sub-sampled dynamic data \eqref{eq:DynPAT} is limited. With data from an experimental phantom that will be described in more detail in Section \ref{sec:ExpPhan}, we were able to show in \cite{ArBeBeCoHuLuOgZh16} that while fbf reconstructions with $M_{{\rm sub}} = 8$ still give acceptable results, using $M_{{\rm sub}} = 16$ leads to reconstructions too heavily impaired by missing-data artefacts and noise. However, an inspection of consecutive frames as shown in Figure \eqref{fig:LimFBF} reveals that the temporal correlation between both noise and artefacts differs strongly from the smooth spatio-temporal evolution of the target. Consequently, noise and artefacts should be effectively removed when using an appropriate smooth spatio-temporal image model. This is the key observation we will utilize to enhance dynamic compressed sensing PAT, either to improve the image quality compared to fbf reconstructions or to allow for higher sub-sampling factors $M_{{\rm sub}}$.

\section{Joint Image Reconstruction and Motion Estimation} \label{sec:JointImRecMotionEst}

\subsection{Simultaneous Motion Estimation} \label{sec:SimMotionEst}

A full non-parametric, spatio-temporal variational scheme reads
\begin{equation}
\hat{p} = \argminsub{p \geqslant0} \Bigg\lbrace  \sum_t^T  \frac{1}{2} \sqnorm{C_t A \, p_t - f_t^c} + \mathcal{R}(p) \Bigg\rbrace \enspace,  \label{eq:SpTempImRec}
\end{equation}
where the regularization $\mathcal{R}(p)$ is now a function of the whole image sequence $p \in \R^{NT}$ that cannot be decomposed over frames, i.e.,  $\mathcal{R}(p) \neq \sum_t \J_t(p_t)$. Here, we choose a particular construction of such a scheme introduced in \cite{BuDiSc16}. For our time-discrete dynamic PAT problem, it is given as:
\begin{equation}
(\hat{p},\hat{v}) = \argminsub{(p \geqslant0, v)} \Bigg\lbrace  \sum_t^T  \frac{1}{2} \sqnorm{C_t A \, p_t - f_t^c} + \alpha \J(p_t) + \beta \mathcal{H}(v_t) + \gamma \mathcal{M}(p,v) \Bigg\rbrace \label{eq:ImRecMoEst}
\end{equation}
Here, each $v_t \in \R^{dN}$ describes a $d$-dimensional vector field describing the motion between $p_t$ and $p_{t+1}$, $\J(p_t)$ and $\mathcal{H}(v_t)$ are spatial regularization terms on image and motion field, respectively, and $\alpha$, $\beta$, $\gamma$ are non-negative regularization parameters. The key term is $\mathcal{M}(p,v)$, which enforces a relation between image sequence $p$ and related motion field sequence $v$ by measuring how well they fulfil a (discretized) motion PDE chosen to model a-priori information about the underlying image dynamics. Note that \eqref{eq:SpTempImRec} can be obtained from \eqref{eq:ImRecMoEst} by dropping $\hat{v}$ from the left hand side and replacing the $\argmin$ over $v$ with a minimization over $v$.

\subsection{Optical Flow Constraints} \label{sec:OptFlow}

The purpose of this work is a proof-of-concept study to show that a more sophisticated spatio-temporal approach like \eqref{eq:ImRecMoEst} can generally improve upon simpler fbf reconstruction. Therefore, we stick to rather generic choices of $\J$, $\mathcal{H}$ and $\mathcal{M}$ and leave the examination of problem-specific regularizers encoding more detailed information about image and dynamics for future work. For $\J$ and $\mathcal{H}$ we choose the popular (isotropic) \termabb{total variation}{TV} functional also used in \cite{ArBeBeCoHuLuOgZh16}. The motion term $\mathcal{M}(p,v)$ should enforce a simple continuity equation, known as the \emph{optical flow} equation \cite{HoSc81} in the field of computer vision:
\begin{equation}
\partial_\tau p(r,\tau) + \left(\nabla_r p(r,\tau) \right) \cdot  v(r,\tau)  = 0 \enspace . \label{eq:OptFlowEq}
\end{equation}
One way to achieve this is to let $\mathcal{M}$ measure the least-squares error of a forward difference discretization of \eqref{eq:OptFlowEq} in time:  
\begin{equation}
\mathcal{M}(p,v) = \sum_t^{T-1} \frac{1}{2} \sqnorm{p_{t+1} - p_t + (\nabla p_t) \cdot v_t} \label{eq:OptFlowReg}
\end{equation}
In total, this leads to the variational scheme
\begin{equation}
\begin{split}
(\hat{p},\hat{v}) = & \argminsub{p \geqslant0, v} \Big\lbrace \E (p,v) \Big\rbrace \mydef \argminsub{p \geqslant0, v} \Bigg\lbrace  \sum_t^T  \frac{1}{2} \sqnorm{C_t A \, p_t - f_t^c} + \\
  &\alpha \|\NabFD p_t \|_1 + \beta \sum_i^d \|\NabFD v_{x_i,t} \|_1 +   \frac{\gamma}{2} \sqnorm{p_{t+1} - p_t + (\NabCD p_t) \cdot v_t} \Bigg\rbrace \enspace, \label{eq:TVTVL2SemiCont}
\end{split}
\end{equation}
where we define $p_{T+1} := p_T$, $v_T := 0$ to simplify the formula. The spatial gradients in the TV terms are implemented with forward differences (denoted by $\NabFD$) as described in Appendix A in \cite{ArBeBeCoHuLuOgZh16}. We chose to implement the TV of the motion field as a sum over the TV of the single components here and leave other possible choices for future work. The spatial gradient $\nabla p_t$ in the optical flow term is discretized using central differences, denoted by $\NabCD$. As the scheme is solved implicitly for $p$ given $v$, this gives a stable discretization of \eqref{eq:OptFlowEq}. More details on the discretization can be found in \cite{Di15}. We will refer to \eqref{eq:TVTVL2SemiCont} as \emph{TVTVL2} model. 


\section{Optimization} \label{sec:Opt}

The TVTVL2 model \eqref{eq:TVTVL2SemiCont} leads to a large-scale, non-smooth, bi-convex optimization problem in $p$ and $v$ involving the computationally intensive acoustic propagation operator $A$ applied to $T$ image frames. We will therefore decompose it into several sub-problems to disentangle its most complicated components. All sub-problems will be solved with iterative first order techniques. For general introductions to numerical optimization suited for imaging applications, we refer to \cite{BuSaSt14,ChPo16}.

\subsection{Forward-Backward Splitting} \label{subsec:FwdBack}

Computing each matrix-vector product $A p_t$ or $A^* f_t$ involves the numerical solution of a potentially inhomogenous 3D wave equation \eqref{eq:PATfwd} with high spatial and temporal resolution. For this, we will use the $k$-space pseudospectral time domain method \cite{MaSoLiTaNaWa01,CoKaArBe07,TrZhCo10} implemented in the \textit{k-Wave} Matlab Toolbox \cite{TrCo10}. With this implementation, each matrix-vector product with $A$ or $A^*$ has the complexity $\mathcal{O}(M_\tau N \log(N))$ \cite{ArBeCoLuTr16} and typically $M_\tau > N$. In contrast, all linear operators in the regularization terms in \eqref{eq:TVTVL2SemiCont} have complexity $\mathcal{O}(N)$. For this reason, we build the outer-most iteration (index $i$) of our scheme by decoupling the smooth, convex data term containing $A$ from all other terms and the non-negativity constraints on $p$ by a \emph{proximal forward-backward splitting}/\emph{proximal gradient descent} scheme (see \cite{GoStBa14} for an extensive overview). For this, we need to define the \emph{proximal operator} of a functional $\J(x)$ as 
\begin{equation}
\prox_{\alpha \J}(y) \mydef \argminsub{x} \left\lbrace \alpha \J(x) + \frac{1}{2} \sqnorm{x - y} \right\rbrace \enspace .\label{eq:Prox}
\end{equation} 
Furthermore, as $v$ is not part of the data term, it appears only in the second step of the iterative scheme:
\begin{subequations}
\begin{align}
\tilde{p}_t &= p_t^i - \eta A^* C_t^* \left(C_t A p^{i}_t - f_t^c  \right) \; \forall\, t = 1,\ldots, T &\textnormal{(forward step)}\phantom{,} \label{eq:ProxGradA}\\
\left(p^{i+1}, v^{i+1} \right) &= \prox_{\eta \mathcal{R}}\left( \tilde{p} \right)   &\textnormal{(backward step)}, \label{eq:ProxGradB}
\end{align}
\end{subequations}
where $\mathcal{R}(p ,v)$ combines all regularization terms on $p$ and $v$ from \eqref{eq:TVTVL2SemiCont}:
\begin{equation}
\mathcal{R}(p, v) :=  \sum_t^T  \alpha \|\NabFD p_t \|_1 + \beta \sum_i^d \|\NabFD v_{x_i,t} \|_1 +   \frac{\gamma}{2} \sqnorm{p_{t+1} - p_t + (\NabCD p_t) \cdot v_t}  \enspace. \label{eq:RegTerms}
\end{equation}
In \eqref{eq:ProxGradA}-\eqref{eq:ProxGradB}, we initialize $p^0 = 0$ and set the step size $\eta$ to $1.5/\max_t L_t$. $L_t$ is an approximation of the Lipschitz constant of $A^* C_t^* C_t A$ which can be pre-computed for a given setting and sub-sampling scheme with a simple power iteration. The basic scheme \eqref{eq:ProxGradA}-\eqref{eq:ProxGradB} is extended by a gradient extrapolation step (\emph{accelerated or fast gradient methods}) which will lead to an asymptotic convergence rate of $\mathcal{O}(1/i^2)$. For this, we use the \emph{FISTA} extrapolation \cite{BeTe09} with restart whenever an increase in the total energy $\E$ occurs. 

\subsection{Biconvex Optimization} \label{sec:BiConv}

Combining \eqref{eq:RegTerms} and \eqref{eq:Prox}, we see that solving the proximal operator in \eqref{eq:ProxGradB} amounts to solving the following \emph{TVTVL2-regularized denoising problem}:
\begin{equation}
\begin{split}
\left(p^{i+1}, v^{i+1} \right) = \prox_{\eta \mathcal{R}}\left( \tilde{p} \right)  = \argminsub{p \geqslant0, v} \Big\lbrace \tilde{\E}(p,v) \Big\rbrace \mydef \argminsub{p \geqslant0, v} \Bigg\lbrace  \sum_t^T  \frac{1}{2} \sqnorm{ p_t - \tilde{p}_t}  \hspace{2em}\\
\hspace{2em} + \eta \alpha \|\NabFD p_t \|_1 + \eta \beta \sum_i^d \|\NabFD v_{x_i,t} \|_1 +   \frac{\eta \gamma}{2} \sqnorm{p_{t+1} - p_t + (\NabCD p_t) \cdot v_t} \Bigg\rbrace \enspace, \label{eq:TVTVL2Denoise}
\end{split}
\end{equation}
The main difficulty here is the motion term. The product $(\NabCD p_t) \cdot v_t$  renders it \emph{bi-convex}, i.e., $\tilde{\E}(p,v)$ is convex in each of the single variables $p$ or $v$ once the other is fixed, but non-convex as a function of both variables. As such, bi-convex problems are global optimization problems that can have a large number of local minima. An overview over bi-convex optimization can be found in \cite{GoPfKl07}. Compared to general global optimization problems, the convex sub-structures can be utilized to design efficient optimization schemes with certain global convergence properties. A popular approach is given by the \termabb{alternate convex search}{ACS} method which alternates between minimizing $\tilde{\E}(p,v)$ for one variable while keeping the other fixed. Applied to \eqref{eq:TVTVL2Denoise}, the ACS iteration (index $j$) reads:
\begin{subequations}
\begin{equation}
p^{j+1} = \argminsub{p \geqslant0} \Bigg\lbrace  \sum_t^T  \frac{1}{2} \sqnorm{ p_t - \tilde{p}_t} + \tilde{\alpha} \|\NabFD p_t \|_1 +   \frac{\tilde{\gamma}}{2} \sqnorm{p_{t+1} - p_t + (\nabla^\pm p_t) \cdot v^{j}_t} \Bigg\rbrace
\label{eq:ACSpOpt}
\end{equation}
\vskip -10pt
\begin{equation}
 v^{j+1} = \argminsub{v} \Bigg\lbrace  \sum_t^T  \tilde{\beta} \sum_i^d \|\NabFD v_{x_i,t} \|_1 +   \frac{\tilde{\gamma}}{2} \sqnorm{p^{j+1}_{t+1} - p^{j+1}_t + (\NabCD p^{j+1}_t) \cdot v_t} \Bigg\rbrace \enspace ,
\label{eq:ACSvOpt}
\end{equation}
\end{subequations}
where we defined $\tilde{\alpha} = \eta \alpha$, $\tilde{\beta} = \eta \beta$, $\tilde{\gamma} = \eta \gamma$. 
The first problem \eqref{eq:ACSpOpt} is a denoising problem for $p$ with a regularization consisting of a TV and a transport term. The problem \eqref{eq:ACSvOpt} for $v$ is an optical flow estimation problem with TV regularization. Note that it is separable in $t$, i.e., it can be solved fbf. Both sub-problems are convex and therefore, approximate solutions can be found reasonably fast by iterative first order schemes (iteration index $k$). The next two sections will present two different approaches for each sub-problem. First, we will repeat how to apply the \termabb{primal dual hybrid gradient}{PDHG} \cite{PoCrBiCh09,ChPo11} algorithm as already proposed in \cite{BuDiSc16}. While this will be sufficient for treating small-scale 2D problems such as examined in Section \ref{sec:NumPhan}, we will then introduce tailored \termabb{alternating direction method of multipliers}{ADMM} (e.g., \cite{BoPaChPeEc11}) schemes that will be shown to be sufficiently efficient to also treat large-scale 3D problems as encountered in the real-data scenarios examined in Section \ref{sec:ExpPhan}. \\
However, as \eqref{eq:ACSpOpt}  and  \eqref{eq:ACSvOpt} are non-smooth, both PDHG and ADMM rely on dual or primal-dual formulations and can therefore not guarantee a monotonous decay of the iterates energy $\tilde{\E}(p,v)$. This leads to a potential problem: While ACS will still converge in objective value $\tilde\E(p,v)$ if we do not solve the sub-problems exactly but only find fast approximate solutions, we need to guarantee that $\tilde{\E}(p,v)$ decreases in every step. Therefore, we will need to track the energies of all iterates and allow sub-routines to run long enough to ensure a sufficient decay. In addition, we will warm-start the sub-routines with all the variables from their last call, even though this will lead to an increased memory consumption.

\subsection{Solution of Convex Subproblems by PDHG} \label{sec:PDHG}

The PDHG algorithm has become the de facto standard template for solving convex, non-smooth optimization problems in a vector space $\mathcal{X}$  involving complicated linear operators $K: \mathcal{X} \rightarrow \mathcal{Y}$ for which matrix-vector products with $K$ and $K^*$ can be computed. The idea is to formulate the problem in the \emph{primal form} as 
\begin{equation}
\min_{x \in \mathcal{X}} \; \E(x) =  \min_{x \in \mathcal{X}} \; \mathcal{G}(x) + \mathcal{F}( K x) \enspace, \label{eq:PDHGprimal}
\end{equation}
with proper, convex functionals $\mathcal{G}$ and $\mathcal{F}$ and to then switch to the equivalent \emph{primal-dual formulation},
\begin{equation}
\min_{x \in \mathcal{X}} \; \max_{y \in \mathcal{Y}} \; \langle Kx , y \rangle + \mathcal{G}(x) - \mathcal{F}^*(y) \enspace, \label{eq:PDHGprimaldual}
\end{equation}
which involves the \emph{convex conjugate} $\mathcal{F}^*$ of $\mathcal{F}$. The advantage of this formulation over \eqref{eq:PDHGprimal} is that the operator $K$ does not show up in the non-linear terms any more. The PDHG algorithm then solves the saddle-point problem \eqref{eq:PDHGprimaldual} by basically alternating a gradient descent in the primal variable and a gradient ascent in the dual variable. In addition, it performs an overrelaxation step in one of the variables (here, the primal one), see Algorithm \ref{algo:PDHG}.
\begin{algorithm} 
\caption{Primal Dual Hybrid Gradient Scheme (PDHG)\label{algo:PDHG}}
Given $\mu > 0, \, \nu > 0, \, \theta \in [0, 1], \, \hat{x}^0, \, y^0$, iterate for $k = 1, 2, \ldots $:
\begin{subequations}
\begin{align}
y^{k+1} &= \prox_{\nu F^*}\left(y^k + \nu K \hat{x}^k \right)   &\textnormal{\small (prox-grad step in $y$)} \label{eq:PDHGa}\\
x^{k+1} &= \prox_{\mu G}\left(x^k - \mu K^* y^{k+1} \right)   & \textnormal{\small (prox-grad step in $x$)}  \label{eq:PDHGb}\\
\hat{x}^{k+1} &= x^{k+1} + \theta \left(x^{k+1} - x^{k}\right)  & \textnormal{\small (overrelaxation)} \label{eq:PDHGc}
\end{align}
\end{subequations}
\end{algorithm}
To apply this to \eqref{eq:ACSpOpt}, i.e., $x = p$, we choose 
\begin{subequations}
\begin{align}
&K p =
\begin{bmatrix}
\NabFD \\
\parFD_t + \sum_i^d v_{x_i}^j \parCD_{x_i} 
\end{bmatrix} 
p \mydef \begin{bmatrix}
\NabFD \\
D_{v^j}
\end{bmatrix} p
\enspace, \qquad 
K^* y = - \NabFD \cdot y_1  + D_{v^j}^* \, y_2
\enspace , \label{eq:K_PDHG_p}\\
&\mathcal{G}(p) = \chi_+(p) + \sum_t^T  \frac{1}{2} \sqnorm{ p_t - \tilde{p}_t} \enspace, \qquad \quad \chi_+(p) \mydef  
\begin{cases}
 0 	&\text{if} \quad p_i \geqslant 0 \; \forall \, i \\
 \infty &\text{else.}
 \end{cases}, \label{eq:G_PDHG_p}\\
&\mathcal{F}(y) = \mathcal{F}\left( 
\begin{bmatrix}
y_1 \\ 
y_2 
\end{bmatrix} 
\right) = \tilde{\alpha} \| y_1 \|_1 +   \frac{\tilde{\gamma}}{2} \sqnorm{y_2} \enspace , \label{eq:F_PDHG_p}
\end{align}
\end{subequations}
where $y_1 \in \R^{N T d}$ represents a $d$-dimensional spatio-temporal vector field resulting from applying $\NabFD$ to every frame of a $d$-dimensional dynamic image $p$ ($d = 2, 3$ here). The explicit form of the proximal operators needed to implement Algorithm \ref{algo:PDHG} with these choices are listed in Appendix \ref{sec:ProxOp}. We can use the PDHG scheme to solve \eqref{eq:ACSvOpt}, i.e., $x = v$, by choosing 
\begin{subequations}
\begin{align}
&K v = \left( I_d \otimes \NabFD \right) 
\begin{bmatrix}
v_{x_1} \\ 
\vdots \\ 
v_{x_d}\\
\end{bmatrix}  \enspace, \qquad \quad K^* y = 
\begin{bmatrix}
- \NabFD \cdot y_{1} \\ 
\vdots \\ 
- \NabFD \cdot y_{d} \\
\end{bmatrix}  \enspace ,  \label{eq:K_PDHG_v}\\
&\mathcal{G}(v) = \frac{\tilde{\gamma}}{2} \sum_t^T \sqnorm{p^{j+1}_{t+1} - p^{j+1}_t + (\NabCD p^{j+1}_t) \cdot v_t} \enspace , \label{eq:G_PDHG_v}\\
&\mathcal{F}(y) = 
\mathcal{F}\left( 
\begin{bmatrix}
y_{1} \\ 
\vdots  \\
y_{d}
\end{bmatrix} 
\right) = \sum_t^T  \tilde{\beta} \sum_i^d \|y_{i,t} \|_1 \enspace . \label{eq:F_PDHG_v}
\end{align}
\end{subequations}
Here, $y \in \R^{N T d^2}$ represents a $(d\times d)$-dimensional spatio-temporal tensor field with components $y_i \in \R^{N T d}$ representing the spatial Jacobian of every frame of a $d$-dimensional dynamic vector field $v$. Again, the solution of the involved proximal operators is shifted to Appendix \ref{sec:ProxOp}. Note that as \eqref{eq:ACSvOpt} can be solved fbf, i.e., for each $v_t$ separately, the PDHG algorithm sketched above can be parallelized over $t$. While this is an appealing option, for its use within ACS one has to implement it such that the overall energy $\tilde{\E}(p,v)$ (which is summed over $t$) deceases sufficiently, cf. Section \ref{sec:BiConv}. \\
The overrelaxation parameter $\theta$ in both PDHG schemes is chosen as $1$. Furthermore, we need to choose the step sizes $\mu$, $\nu$ in dependence on $K$ to ensure convergence (cf. \cite{ChPo11,ChPo16}). Due to the complicated structure of \eqref{eq:K_PDHG_p}, we use the extension of PDHG by diagonal preconditioning proposed in \cite{PoCh11} (the $\alpha$ parameter in \cite{PoCh11} is set to $1$), which is easy to compute for our problem and was found to work well compared to standard choices based on estimates of $\norm{K}_{2,2}$. The operator $K$ in \eqref{eq:K_PDHG_v} has a simple structure and it can be shown that $\norm{K}_{2,2} \leqslant 4d$ \cite{ChPo11}. As such, the choice $\mu = 1/2d$, $\nu = 1/2$ fulfils $\mu \nu \norm{K}_{2,2}^2 \leqslant 1$ and leads to convergence (this balancing between $\mu$ and $\nu$ was found empirically).

\subsection{Solution of Convex Subproblems by ADMM} \label{sec:ADMM}

In the ADMM approach, the unconstrained but coupled convex problem \eqref{eq:PDHGprimal} is first converted into an equality-constrained but uncoupled convex problem by introducing an auxiliary variable $y = Kx$,
\begin{equation}
\eqref{eq:PDHGprimal} \Longleftrightarrow \min_{x \in \mathcal{X} , y \in \mathcal{Y}} \; \mathcal{G}(x) + \mathcal{F}(y) \quad such \; that \quad y = Kx \enspace, \label{eq:ADMMsplit}
\end{equation}
which is then solved by a combination of \emph{dual ascent}, \emph{augmented Lagrangian techniques}, and the \emph{method of multipliers}. The final ADMM scheme is described in Algorithm \ref{algo:ADMM}.
\begin{algorithm} 
\caption{Alternating Direction Method of Multipliers (ADMM)\label{algo:ADMM}}
Given $\rho > 0$, \, $y^0$, \, $w^0$, iterate for $k = 1, 2,\ldots$:
\begin{align}
x^{k+1} &= \argminsub{x \in \mathcal{X}} \left\lbrace \mathcal{G}(x) + \frac{\rho}{2} \sqnorm{K x - y^k + w^k} \right\rbrace \label{eq:ADMMx}\\
y^{k+1} &= \argminsub{y \in \mathcal{Y}} \left\lbrace \mathcal{F}(y) + \frac{\rho}{2} \sqnorm{K x^{k+1} - y + w^k} \right\rbrace \label{eq:ADMMy}\\
w^{k+1} &= w^k + K x^{k+1} - y^{k+1} \label{eq:ADMMw}
\end{align}
\end{algorithm}
The crucial difference to the PDHG schemes is that the update of $x$, \eqref{eq:ADMMx}, is now implicit, and we will choose the split $y = Kx$ such that it will be given as the solution of a least-squares problem involving all linear operators\footnote{ADMM and PDHG schemes are actually very closely related, although most introductions of the two methods do not immediately imply this, and our short overview here cannot cover it. See  \cite{BuSaSt14,ChPo16} for an extensive discussion.}. This can be advantageous in cases where $K$ suffers from bad conditioning, but only leads to a computationally efficient scheme if the corresponding normal equations can be solved fast. Fortunately, ADMM still converges if the sub-problems \eqref{eq:ADMMx} and \eqref{eq:ADMMy} are solved approximately but with accuracy increasing with $k$ (see \cite{Es09} and references therein for a precise statement). Therefore, warm-started iterative linear solvers with carefully chosen stop conditions can be used. For problem \eqref{eq:ACSpOpt}, i.e., $x = p$, we realize the ADMM iteration by  
\begin{subequations}
\begin{align}
 K p &=
\begin{bmatrix}
\NabFD \\ 
I_N 
\end{bmatrix}  p \qquad \Rightarrow \qquad
K^* 
\begin{bmatrix}
y_1 \\ 
y_2 
\end{bmatrix} 
 = - \NabFD \cdot y_1  + y_2
\enspace , \label{eq:K_ADMM_p} \\
\mathcal{G}(p) &= \sum_t^T  \frac{1}{2} \sqnorm{ p_t - \tilde{p}_t} + \frac{\tilde{\gamma}}{2} \sqnorm{p_{t+1} - p_t +  v^{j}_t \cdot (\nabla p_t)} \label{eq:G_ADMM_p}  \\
&= \frac{1}{2} \sqnorm{
\begin{bmatrix}
I_N \\
\sqrt{\tilde{\gamma}} D_{v^j}
\end{bmatrix} p - 
\begin{bmatrix}
\tilde{p} \\
0
\end{bmatrix} 
}  \enspace , \nonumber\\
\mathcal{F}(y) &= \mathcal{F}\left( 
\begin{bmatrix}
y_1 \\ 
y_2
\end{bmatrix} 
\right) = \tilde{\alpha} \| y_1 \|_1 + \chi_+(y_2) \enspace , \label{eq:F_ADMM_p}
\end{align}
\end{subequations}
where $y_1 \in \R^{N T d}$ represents a $d$-dimensional\ spatio-temporal vector field resulting from applying $\NabFD$ to a $d$-dimensional\ dynamic image $p$ and $y_2 \in \R^{N T}$ accounts for the non-negativity constraints. We will denote the corresponding parts of $w$ by $w_1$ and $w_2$ as well. For these choices, the update \eqref{eq:ADMMx} is given by
\begin{multline}
p^{k+1} = \argminsub{p} \left\lbrace \frac{1}{2} \sqnorm{
\begin{bmatrix}
I_N \\
\sqrt{\tilde{\gamma}} \; D_{v^j}\\
\sqrt{\rho} \; \NabFD \\ 
\sqrt{\rho} \; I_N 
\end{bmatrix}  p - 
\begin{bmatrix}
\tilde{p} \\
0 \\
\sqrt{\rho} \left( y^k_1 - w^k_1 \right)\\ 
\sqrt{\rho} \left( y^k_2 - w^k_2 \right)
\end{bmatrix}
} \right\rbrace  \label{eq:LinSys_p}\\
 = \left( (1 + \rho) I_N + \tilde{\gamma} D_{v^j}^* D_{v^j} + \rho \laplaceFD \right)^{-1}  \left(\tilde{p} + \rho \NabFD \cdot \left( y^k_1 - w^k_1 \right) + \rho \left( y^k_2 - w^k_2 \right) \right) \enspace .
\end{multline}
All the linear operators can easily be implemented in a matrix-free way and so \eqref{eq:LinSys_p} can be solved with a standard \termabb{conjugate gradient}{CG} implementation. As in the PDHG schemes, update \eqref{eq:ADMMy} can be solved explicitly using proximal operators (Appendix \ref{sec:ProxOp}). \\
For solving \eqref{eq:ACSvOpt}, i.e., $x = v$, choose exactly the same split as in the corresponding PDHG scheme, i.e., \eqref{eq:K_PDHG_v}-\eqref{eq:F_PDHG_v} and make use of the fact that the optimization can be solved for each $t$ separately: For each $t$, the update \eqref{eq:ADMMx} is given by
\begin{equation}
v^{k+1}_t = \left( \tilde{\gamma} E^* E + \rho I_d \otimes  \laplaceFD  \right)^{-1} \left(  \tilde{\gamma} E^* \left(- p^{j+1}_{t+1} + p^{j+1}_t\right) + \rho K^* \left(y^k - w^k \right) \right) \enspace , \label{eq:LinSys_v}
\end{equation}
where $E = E(p^{j+1}_t)$ is an $N \times dN$ matrix implementing the point-wise multiplication and summation of the components of a vector field with the spatial gradients of $p^{j+1}_t$, i.e.,  $E v = \sum_i^d v_{x_i} \parCD_{x_i} p^{j+1}_t$ as
\begin{equation}
E = 
\begin{bmatrix}
\diag \left( \parCD_{x_1} p^{j+1}_t \right) & \ldots & \diag \left( \parCD_{x_d} p_t \right) 
\end{bmatrix}  \enspace .
\end{equation}
For $d=2$, we take a closer look at the structure of the matrix to invert in \eqref{eq:LinSys_v}:
\begin{equation}
\left( \tilde{\gamma} E^* E + \rho I_d \otimes  \laplaceFD  \right) = \tilde{\gamma} 
\begin{bmatrix}
T_{11} & T_{12} \\
T_{11} & T_{22} 
\end{bmatrix} +  \rho 
\begin{bmatrix}
\laplaceFD & 0 \\
0 & \laplaceFD
\end{bmatrix} \enspace ,
\end{equation}
where $T_{kl} \mydef \diag \left( \parCD_{x_k} p^{j+1}_t \cdot \parCD_{x_l} p^{j+1}_t \right)$. While one can easily implement matrix-free iterative solvers for this system, we chose to explicitly build this very sparse matrix to be able to use efficient pre-conditioning techniques. Within the ADMM iteration, this comes with little overhead as only the right hand side in system \eqref{eq:LinSys_v} changes during the iteration. We will examine different combinations of pre-conditioners and iterative solvers in the numerical studies \cite{Sa03}:
As pre-conditioners, we consider
\begin{itemize}
\item IC(0): incomplete Cholesky pre-conditioner with zero-fill as implemented in \textsc{Matlab} (R2016a). 
\item ICT: incomplete Cholesky pre-conditioner with threshold dropping (threshold: 1e-3) as implemented in \textsc{Matlab} (R2016a). 
\item AMG: Algebraic Multigrid W-cycle pre-conditioner based on the implementation in \cite{Ch09}, which uses a modification of Ruge-Stuben coarsening, two-points interpolation (use at most two connected coarse nodes) and a  direct solver on the coarsest level. 
\end{itemize}
As iterative solvers, we examine the standard CG method and the \termabb{Minimum Residual Method}{MINRES}. Further details will be discussed in the  next section. The iterative solvers are warm-started with the previous solution $p^k$ or $v_t^k$, perform at least 3 iterations and stop when the relative residual norm is below $tol(k) = 10^{-3}/k^{3/2}$, i.e., we progressively increase the precision to which we solve sub-problem \eqref{eq:ADMMx}. Update \eqref{eq:ADMMy} can be solved as for the corresponding PDHG scheme (Appendix \ref{sec:ProxOp}).\\
While ADMM converges for all $\rho > 0$, its choice has a crucial impact on the speed of convergence and other properties of the iterates, e.g., the monotonicity of the energy $\mathcal{E}\left(x^{k}\right)$, which is important for using ADMM inside of an ASC (cf. Section \ref{sec:BiConv}). For using ADMM on the $p$ update \eqref{eq:ACSpOpt}, we use the adaptation strategy described in Section 3.4.1 of \cite{BoPaChPeEc11} during the steps $k = 1,\ldots,25$ and fix it thereafter. For the first $p$ update \eqref{eq:ACSpOpt} within ASC we initialize $\rho = 1$ and then always warm-start the following $p$ update with the adapted $\rho$. In the ADMM scheme for the $v$ update \eqref{eq:ACSvOpt}, we fix $\rho = 10^{-1}$ for $d=2$ and $10^2$ for $d=3$, firstly to avoid a re-computation of the matrices and their pre-conditioners (see above) and secondly to enforce a fast transition to the regime of monotonous energy decay. As with any alternating optimization, the ADMM scheme can benefit from over-relaxation. We use the technique discussed in Section 3.4.3 of \cite{BoPaChPeEc11}, which  consists of replacing the quantity $K x^{k+1}$ in Algorithm \ref{algo:ADMM} by $s Kx^{k+1} + (1-s) y^k$. Throughout the experiments, we use $s = 1.8$.\\
Although we limited our presentation here to the most important features, it already became apparent that compared to PDHG, ADMM schemes are more difficult to design and parametrize. Also note that ADMM with the specific type of split that we used here is equivalent to the \emph{split Bregman method} \cite{GoOs09,Es09}, which derives Algorithm \ref{algo:ADMM} from a different perspective.


\section{Results} \label{sec:Res}

In this section, we first demonstrate the main features of the proposed methods on a simple numerical phantom in 2D before we discuss their realization for experimental data in 3D. As we can only show snapshots for a few time frames of the reconstructions here, movies of all reconstructions can be found in the supplementary material. For computing the results presented, we used ADMM in both the $p$ update \eqref{eq:ACSpOpt} and the $v$ update \eqref{eq:ACSvOpt} as described in the previous section. In the $v$ update, AMG-CG was used as a least squares solver. In Section \ref{sec:ExpOpt}, we compare this choice to possible alternatives in more detail. \\
All routines have been implemented as part of a Matlab toolbox for PAT image reconstruction which will be made available in near future. The toolbox relies on the k-Wave toolbox (see \cite{TrCo10},  \href{http://www.k-wave.org/}{http://www.k-wave.org/}) to implement $A$ and $A^*$, which allows to use highly optimized C++ and CUDA code to compute the 3D wave propagation on parallel CPU or GPU architectures.

\subsection{Numerical 2D Phantom} \label{sec:NumPhan}

\begin{figure}[tb]
   \centering
\hfill
\subfloat[][\label{subfig:DynEll_13}]{\includegraphics[height=0.4\textwidth]{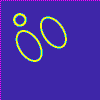}
\hskip 0.01\textwidth
\fbox{\includegraphics[height=0.4\textwidth]{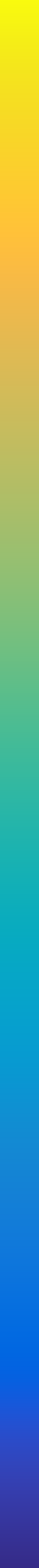}}}
\hfill
\subfloat[][\label{subfig:DynEll_color}]{\fbox{\includegraphics[height=0.4\textwidth]{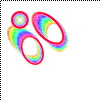}}
\hskip 0.01\textwidth
\includegraphics[height=0.4\textwidth]{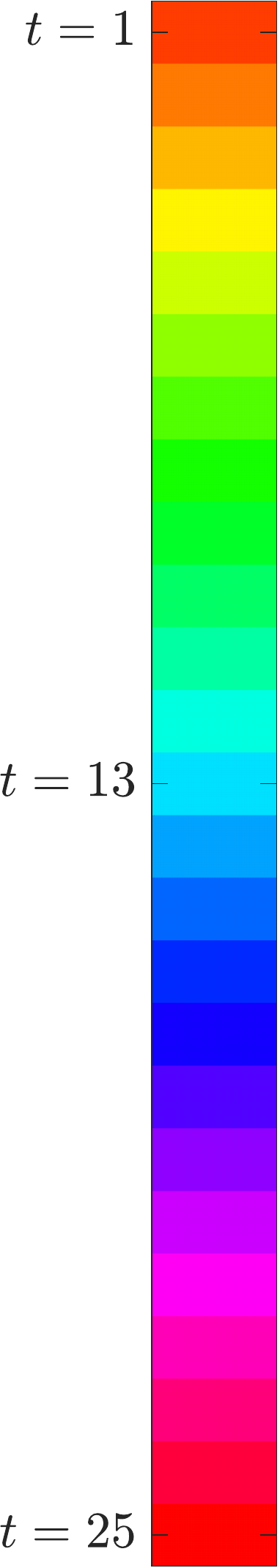}}
\caption{The 2D numerical phantom. \protect\subref{subfig:DynEll_13}: A snapshot at $t = 13$ and the corresponding color bar. The sensor locations are shown as pink pixels (left and top edge). \protect\subref{subfig:DynEll_color} A visualization of all $T = 25$ frames as a color-coded RGB overlay. The color bar displays which colors are assigned to which time frame, the sensor pixels are shown as black pixels.}
   \label{fig:DynEll_setup}
\end{figure}

\begin{figure}[tb]
   \centering
\subfloat[][phantom $p$ (ground truth)]{\includegraphics[width=0.32\textwidth]{DynamicEllipses_t13.png}}
\hfill
\subfloat[][NNLS, cnv \label{subfig:DynEll_LS_13}]{\includegraphics[width=0.32\textwidth]{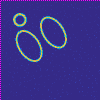}}
\hfill
\subfloat[][TV-fbf, $\alpha = 2 \cdot 10^{-3}$ cvn\label{subfig:DynEll_TV_13}]{\includegraphics[width=0.32\textwidth]{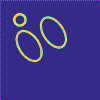}} \\
\subfloat[][]{\includegraphics[width=0.32\textwidth]{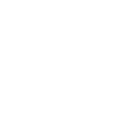}}
\hfill
\subfloat[][NNLS, rSP-25\label{subfig:DynEll_LSSS_13}]{\includegraphics[width=0.32\textwidth]{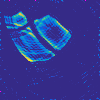}}
\hfill
\subfloat[][TV-fbf, $\alpha = 3.2 \cdot 10^{-4}$ rSP-25\label{subfig:DynEll_TVSS_13}]{\includegraphics[width=0.32\textwidth]{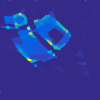}} \\
\caption{Snapshots at $t = 13$ of the results of frame-by-frame image reconstruction methods \eqref{eq:ImRecFBF} for full (cvn) and sub-sampled (rSP-25) data.}
   \label{fig:DynEll_fbf}
\end{figure}

\begin{figure}[tb]
   \centering
\subfloat[][phantom $p$ (ground truth)]{\includegraphics[width=0.32\textwidth]{DynamicEllipses_t13.png}}
\hfill
\subfloat[][TVTVL2, $\gamma = 1$, cnv \label{subfig:DynEll_TVTVL2_para1_13}]{\includegraphics[width=0.32\textwidth]{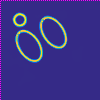}}
\hfill
\subfloat[][TVTVL2, $\gamma = 0.1$, cnv\label{subfig:DynEll_TVTVL2_para2_13}]{\includegraphics[width=0.32\textwidth]{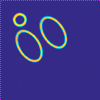}} \\
\subfloat[][]{\includegraphics[width=0.32\textwidth]{dummy.png}}
\hfill
\subfloat[][TVTVL2, $\gamma = 1$, rSP25  \label{subfig:DynEll_TVTVL2SS_para1_13}]{\includegraphics[width=0.32\textwidth]{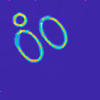}}
\hfill
\subfloat[][TVTVL2, $\gamma = 0.1$, rSP-25\label{subfig:DynEll_TVTVL2SS_para2_13}]{\includegraphics[width=0.32\textwidth]{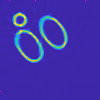}} 
\caption{Snapshots at $t = 13$ of the results of the TVTVL2 image reconstruction \eqref{eq:TVTVL2SemiCont} for full (cvn) and sub-sampled (rSP-25) data. The parameters $\alpha$ and $\beta$ were set to the corresponding value of the $\alpha$ used for the TV-fbf reconstructions in  Figure \ref{fig:DynEll_fbf} and $\gamma$ was set to $1$ or $0.1$.}
   \label{fig:DynEll_TVTVL2}
\end{figure}

\begin{figure}[tb]
   \centering
\subfloat[][$(\hat{\alpha},\hat{\alpha}, 1)$]{\includegraphics[width=0.24\textwidth]{DyEl_TVTVL2_a32e-05b32e-05g1e+00_t13_SS.png}}
\hfill
\subfloat[][$(\hat{\alpha},\hat{\alpha}/4, 1)$]{\includegraphics[width=0.24\textwidth]{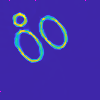}}
\hfill
\subfloat[][$(\hat{\alpha}/4,\hat{\alpha}/4, 1)$]{\includegraphics[width=0.24\textwidth]{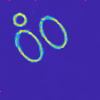}}
\hfill
\subfloat[][$(\hat{\alpha}/4,\hat{\alpha}/16, 1)$]{\includegraphics[width=0.24\textwidth]{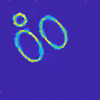}}\\
\subfloat[][$(\hat{\alpha},\hat{\alpha}, 0.1)$]{\includegraphics[width=0.24\textwidth]{DyEl_TVTVL2_a32e-05b32e-05g1e-01_t13_SS.png}}
\hfill
\subfloat[][$(\hat{\alpha},\hat{\alpha}/4, 0.1)$]{\includegraphics[width=0.24\textwidth]{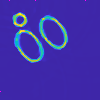}}
\hfill
\subfloat[][$(\hat{\alpha}/4,\hat{\alpha}/4, 0.1)$]{\includegraphics[width=0.24\textwidth]{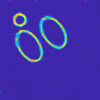}}
\hfill
\subfloat[][$(\hat{\alpha}/4,\hat{\alpha}/16, 0.1)$]{\includegraphics[width=0.24\textwidth]{DyEl_TVTVL2_a8e-05b8e-05g1e-01_t13_SS.png}}
\caption{Snapshots at $t = 13$ of the results of the TVTVL2 image reconstruction \eqref{eq:TVTVL2SemiCont} for sub-sampled (rSP-25) data for different combinations of $(\alpha, \beta, \gamma)$. Here, $\hat \alpha$ corresponds to the value of the regularization parameter used for the TV-fbf reconstructions in  Figure \ref{fig:DynEll_fbf}.}
   \label{fig:DynEll_TVTVL2_ParaCmp}
\end{figure}

\begin{figure}[tb]
   \centering
\subfloat[][phantom $p$ (reference) \label{subfig:DynEll_v_13}]{\includegraphics[width=0.32\textwidth]{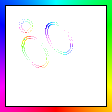}}
\hfill
\subfloat[][TVTVL2, $\gamma = 1$, cnv \label{subfig:DynEll_v_TVTVL2_para1_13}]{\includegraphics[width=0.32\textwidth]{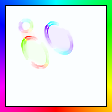}}
\hfill
\subfloat[][TVTVL2, $\gamma = 0.1$, cnv\label{subfig:DynEll_v_TVTVL2_para2_13}]{\includegraphics[width=0.32\textwidth]{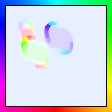}}\\ 
\subfloat[][color scheme for 2D vectors \label{subfig:ColourWheel}]{\includegraphics[width=0.32\textwidth]{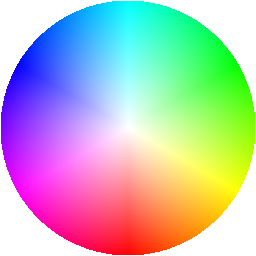}}
\hfill
\subfloat[][TVTVL2, $\gamma = 1$, rSP25 \label{subfig:DynEll_v_TVTVL2SS_para1_13}]{\includegraphics[width=0.32\textwidth]{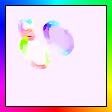}}
\hfill
\subfloat[][TVTVL2, $\gamma = 0.1$, rSP25\label{subfig:DynEll_v_TVTVL2SS_para2_13}]{\includegraphics[width=0.32\textwidth]{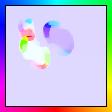}} \\
\caption{Illustration of the reconstructed motion fields $v$ for the results shown in Figure \ref{fig:DynEll_TVTVL2}. First, each vector field is rescaled such that $\max_i \norm{v_i}_2  = 1$. Then, each pixel $i$ is coloured by mapping  direction and norm of $v_i$ to the color scheme displayed in \protect\subref{subfig:ColourWheel}. To further ease the visualization, a coloured frame was added to each motion image to depict the colors corresponding to a vector pointing from the middle of the  image to the pixels of the frame. The motion field for the phantom \protect\subref{subfig:DynEll_v_13} was computed by solving \eqref{eq:ACSvOpt} with the true $p$ as input and $\gamma = 1$, $\beta = 10^{-6}$. }
   \label{fig:DynEll_v_TVTVL2}
\end{figure}

The computational domain is a square of length $20$\si{\milli \meter} which is divided into $N = 100 \times 100$ pixels. Its acoustic properties are assumed homogeneous with $c = 1500$ \si{\meter\per\second}. The conventionally scanned (fully sampled) measurement data (referred to as ``cnv'') is acquired at $M = 100$ sensors sampled at $M_\tau = 472$ time steps with $\delta_\tau = 40$ \si{\nano\second}. The sensors are arranged in two orthogonal lines, which corresponds to a 2D version of a scanning system using two orthogonal Fabry-P\'erot sensors \cite{ElOgZhBeCo16,FoMaLuElAnArBeCo17}. This way, reconstructions from the fully sampled sensor array will not suffer from severe limited view artefacts and we can concentrate on the effects of sub-sampling. \\
The dynamical phantom consists of three tubes that change center position, orientation and size smoothly over $T= 25$ frames and should loosely resemble the dynamics of the X-slices of the experimental phantom in Figure \ref{fig:LimFBF}. Figure \ref{fig:DynEll_setup} shows different visualizations of the phantom. White Gaussian noise with a standard deviation of $\sigma = 5\cdot 10^{-3}$ was added to the simulated pressure time series leading to an average SNR of 20.65 \si{\decibel}. To sub-sample the data, we now assume that in each frame, we can acquire data at a sub-set of $4$ out of the $100$ original sensor locations that have been chosen random but disjoint, such that after $T = 25$ frames, each location has been scanned once. This means that the sub-sampling factor $M_{{\rm sub}} = M/M_c$ (cf. Section \ref{subsec:DynPAT}) equals $T$, i.e., we acquire all $25$ frames with the same scanning time as a single frame in the full data set-up. The operators $C_t$ can thus be written as $C_t = I_{M_\tau} \otimes \tilde{C}_t$ with $\tilde{C}_t$ being a binary $4 \times 100$ matrix with $4$ ones on the main diagonal and all zero otherwise and
\begin{equation}
\bar{C} = 
\begin{bmatrix}
C_{1} \\
\vdots \\
C_{T}
\end{bmatrix}  \label{eq:rSP}
\end{equation}
is a row-permutation of $I_{M}$. We will denote this sub-sampling strategy by rSP-25.\\
First, we compute fbf reconstructions \eqref{eq:ImRecFBF} without any regularizer $\J(p)$, i.e., \termabb{non-negative least squares}{NNLS}, and then using a TV functional (denoted as \emph{TV-fbf}). For this, we use $100$ iterations of the accelerated proximal gradient descend introduced in Section \ref{subsec:FwdBack}. The proximal step \eqref{eq:ProxGradB} is simply a projection onto the positive orthant for NNLS, while it amounts to solving a TV-regularized denoising problem in the case of TV (for details, see \cite{ArBeBeCoHuLuOgZh16}). The results are shown in Figure \ref{fig:DynEll_fbf} and again demonstrate that while we can obtain a good reconstruction with fbf methods for full data, they fail for severely sub-sampled data, similar to the motivating example shown in Figure \ref{fig:LimFBF}. Next we compute reconstructions with the TVTVL2 model \eqref{eq:TVTVL2SemiCont}: An apparent challenge of this more sophisticated spatio-temporal model that we did not discuss up to now is that it relies on three regularization parameters $\alpha$, $\beta$ and $\gamma$. For TV-fbf, it is easy to fix the single parameter $\alpha$ manually: We computed reconstructions for different $\alpha$ for a single frame, and then used the smallest $\alpha$ that visually removed most noise for all frames, which we will denote as $\hat{\alpha}$. For the TVTVL2 model, we start with simply setting $\alpha = \beta = \hat{\alpha}$ and $\gamma = 1$. Figures \ref{subfig:DynEll_TVTVL2_para1_13} and \ref{subfig:DynEll_TVTVL2SS_para1_13} show the results of this naive parameter choice. Although the reconstructions for the sub-sampled data still suffer from some blurring and artefacts, one can clearly see a significant improvement compared to the fbf reconstructions in Figure \ref{fig:DynEll_fbf}. We then varied $(\alpha, \beta, \gamma)$ around this first guess. Figures \ref{subfig:DynEll_TVTVL2_para2_13} and \ref{subfig:DynEll_TVTVL2SS_para2_13} show the effect of decreasing $\gamma$ to $0.1$ which has by far the biggest positive impact. Figure \ref{fig:DynEll_TVTVL2_ParaCmp} illustrate the effects of also varying $\alpha$ and $\beta$, which leads to trade-offs between over-smoothing and artefact reduction. We leave a more detailed parameter study for future work and instead investigate the estimated motion fields. Figure \ref{fig:DynEll_v_TVTVL2} shows that main features of the motion fields can be re-constructed even from sub-sampled data. In particular, the motion fields facilitate the distinction and tracking of different moving objects.

\subsection{Experimental 3D Phantom} \label{sec:ExpPhan}

Now, we examine the performance of the methods on high resolution 3D reconstructions from dynamic experimental phantom data. As outlined in the motivation in Section \ref{subsec:PreWork}, we use the same data as in \cite{ArBeBeCoHuLuOgZh16}, to investigate if the methods described here can improve upon the fbf reconstructions for $M_{{\rm sub}} = 16$ (cf. Figure \ref{fig:LimFBF}). However, for this article to be self-contained, we first briefly recap the set-up and pre-processing used.

\subsubsection{Setup and Pre-processing} \label{sec:ExpSetup}

\begin{figure}[tb]
   \centering
\includegraphics[width = \textwidth]{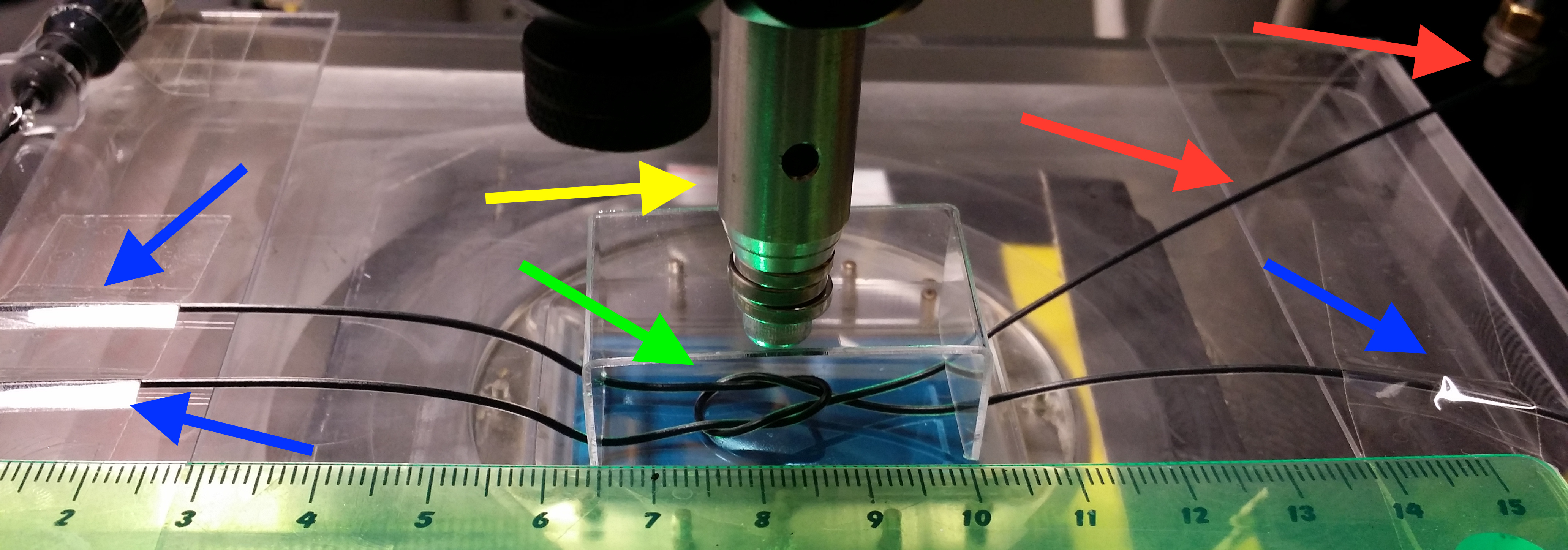}
\caption{Experimental setup for dynamic, stop-motion phantom: The two polythene tubes are immersed in a bath of intralipid solution placed on the FP sensor plane (green arrow). Three of the tube ends are fixated (blue arrows) while one is tied to a motor shaft (red arrows). The excitation laser (yellow arrow) is illuminating from the top.  \label{fig:ThKnSetup}}
   \label{fig:ExpPhantoms}
\end{figure}

The phantom consists of two polythene tubes filled with 100\% and 10\% ink immersed in a $1$\% Intralipid solution with de-ionised water. The tubes were interleaved to form a knot with 4 open ends. As shown in Figure \ref{fig:ThKnSetup}, while three of the ends are fixated, one is tied to a motor shaft. We then acquired PA data using a FP scanner in a stop-motion style: With the whole arrangement at rest, a full, conventional scan was performed. Then, the motor shaft was turned by a fixed angle which caused the knot to both move towards the motor and tighten, and the new arrangement is scanned again. In total, $T = 45$ frames were acquired. The excitation laser pulses were delivered at a rate of $20$\si{\hertz}, had a wavelength of $1064$\si{\nano\meter} and an energy of around $20$\si{\milli\joule}. For a full, conventional scan, pressure time courses at $134 \times 133$ locations on a spatial grid with grid size $150$\si{\micro\meter} were measured for $M_t = 625$ time points with a temporal resolution of $12$\si{\nano\second}. For preprocessing, the data was first clipped to $132 \times 132$ locations. Then, we preformed baseline-correction, band-pass filtering ($0.5$-$20$\si{\mega\hertz}), noisy-channel exclusion and clipped the time courses to the time points $10-400$. More details can be found in \cite{ArBeBeCoHuLuOgZh16}. Note that the signal recorded by the FP sensor is only proportional to the acoustic pressure. To obtain absolute pressure values, one would need to calibrate it with an ultrasound transducer prior to the measurement. While this is necessary to perform quantitative, spectroscopic inference in a second analysis step \cite{CoLaArBe12,FoMaLuElAnArBeCo17}, we did not do it here and all images shown can be considered in arbitrary units. \\
For the inversion, we assume a homogenous sound speed of $1540$ \si{\meter \per \second} and use a 3D spatial grid of dimensions $44 \times 264 \times 264$  with grid size $75$\si{\micro\meter} (the reason for this up-sampling in space is the over-sampling in time and explained in \cite{ArBeBeCoHuLuOgZh16}). Reconstructions from the full, conventional data will again be denoted by "cnv" and will be used to provide a ground truth. The sub-sampled data is generated using the same scheme  \eqref{eq:rSP} as for the simulated data, except that $M_{{\rm sub}} = 16$. The sub-sampling operators are repeated periodically, i.e., $C_{M_{{\rm sub}}+i} = C_{i}$.

\subsubsection{Experimental Results} \label{sec:ExpRes}

\begin{figure}[tb]
   \centering
\subfloat[][TV cnv, X mIP \label{subfig:ExpTVX}]{\includegraphics[width=0.48\textwidth]{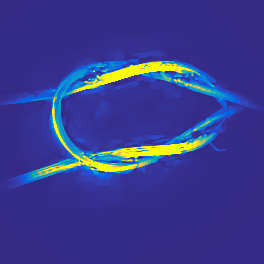}}
\hfill
\subfloat[][TVTVL2 cnv, X mIP \label{subfig:ExpTVTVL2X}]{\includegraphics[width=0.48\textwidth]{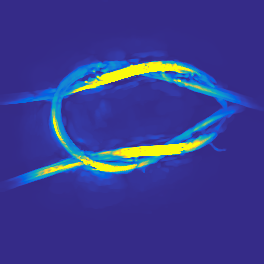}}\\[-8pt]
\subfloat[][TV cnv, Y mIP \label{subfig:ExpTVY}]{\includegraphics[width=0.48\textwidth]{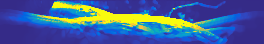}}
\hfill
\subfloat[][TVTVL2 cnv, Y mIP \label{subfig:ExpTVTVL2Y}]{\includegraphics[width=0.48\textwidth]{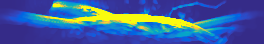}}\\[-8pt]
\subfloat[][TV cnv, Z mIP \label{subfig:ExpTVZ}]{\includegraphics[width=0.48\textwidth]{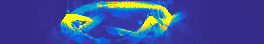}}
\hfill
\subfloat[][TVTVL2 cnv, Z mIP \label{subfig:ExpTVTVL2Z}]{\includegraphics[width=0.48\textwidth]{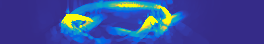}}\\[-8pt]
\subfloat[][TV rSP16, X mIP \label{subfig:ExpTVSSX}]{\includegraphics[width=0.48\textwidth]{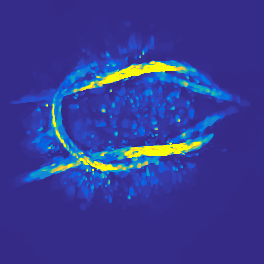}}
\hfill
\subfloat[][TVTVL2 rSP16, X mIP \label{subfig:ExpTVTVL2SSX}]{\includegraphics[width=0.48\textwidth]{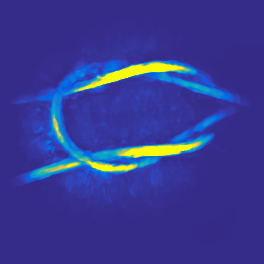}}\\[-8pt]
\subfloat[][TV rSP16, Y mIP \label{subfig:ExpTVSSY}]{\includegraphics[width=0.48\textwidth]{TV-FbF-SS_t23_mIY.png}}
\hfill
\subfloat[][TVTVL2 rSP16, Y mIP \label{subfig:ExpTVTVL2SSY}]{\includegraphics[width=0.48\textwidth]{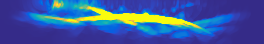}}\\[-8pt]
\subfloat[][TV rSP16, Z mIP \label{subfig:ExpTVSSZ}]{\includegraphics[width=0.48\textwidth]{TV-FbF-SS_t23_mIZ.png}}
\hfill
\subfloat[][TVTVL2 rSP16, Z mIP \label{subfig:ExpTVTVL2SSZ}]{\includegraphics[width=0.48\textwidth]{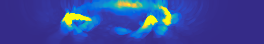}}
   \caption{Snapshots at $t = 23$ of the reconstructed pressure $p$ for full (cvn) and sub-sampled  (rSP-16)  experimental data. \label{fig:Exp_mIP}}
\end{figure}

\begin{figure}[tb]
   \centering
\subfloat[][cnv, slice $x = 19$ \label{subfig:ExpVX}]{\includegraphics[width=0.48\textwidth]{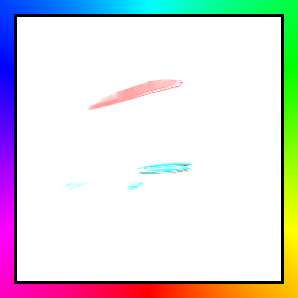}}
\hfill
\subfloat[][rSP-16, slice $x = 19$ \label{subfig:ExpSSVX}]{\includegraphics[width=0.48\textwidth]{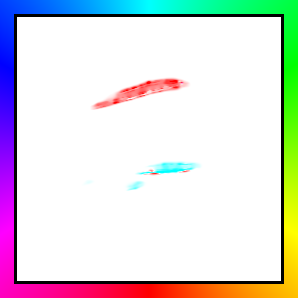}}\\
%
%
\subfloat[][cnv, slice $z = 132$ \label{subfig:ExpVZ}]{\includegraphics[width=0.48\textwidth]{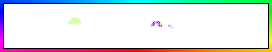}}
\hfill
\subfloat[][rSP-16, slice $z = 132$ \label{subfig:ExpSSVZ}]{\includegraphics[width=0.48\textwidth]{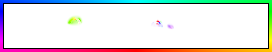}}
   \caption{Snapshots of the non-translational part of the motion field $v^{\rm NT} = v - N^{-1} \sum_i v_i$ reconstructed by the TVTVL2 method for full (cvn) and sub-sampled  (rSP-16)  experimental data between frames $t=22$ and $t=23$. For each image, a single slice along a particular dimension was extracted and only the components of the vector field in the remaining two dimensions is depicted here. The resulting 2D vector field is color-coded in the same way as in Figure \ref{fig:DynEll_v_TVTVL2}. \label{fig:Exp_v_slices}}
\end{figure}

We used the same strategy to choose the regularization parameters as before: For the TVTVL2 model, we choose $\alpha = \beta =  \hat{\alpha}$, $\gamma = 0.1$, where $\hat{\alpha}$ is the regularization parameter for TV-fbf that yields a good compromise between removing noise, sub-sampling and image features (cf. Figure \ref{fig:LimFBF}). Figure \ref{fig:Exp_mIP} shows the results after $20$ iterations (index $i$) of the accelerated proximal gradient descend. Again, we can see a significant improvement of using the simultaneous motion estimation introduced by TVTVL2 compared to TV-fbf. The motion of our phantom has two dominant components: a translation component resulting from pulling the whole knot towards the motor shaft by one tube end, and a component describing the contraction resulting from the three other tube ends being fixed. To examine the later component, we suppress the translation by subtracting the mean motion vector in every frame $\bar{v}_t = N^{-1} \sum_i {(v_t)}_i$. Figure \ref{fig:Exp_v_slices} shows the remaining parts of the motion. Both the fields reconstructed from full and from sub-sampled data accurately describe the contraction. The coloring indicates that the tubes move towards each other, i.e., the knot contracts. 

\subsection{Optimization} \label{sec:ExpOpt}

\begin{figure}[tb]
   \centering
\subfloat[][2D scenario\label{subfig:Cmp2D}]{\includegraphics[height=0.46\textwidth]{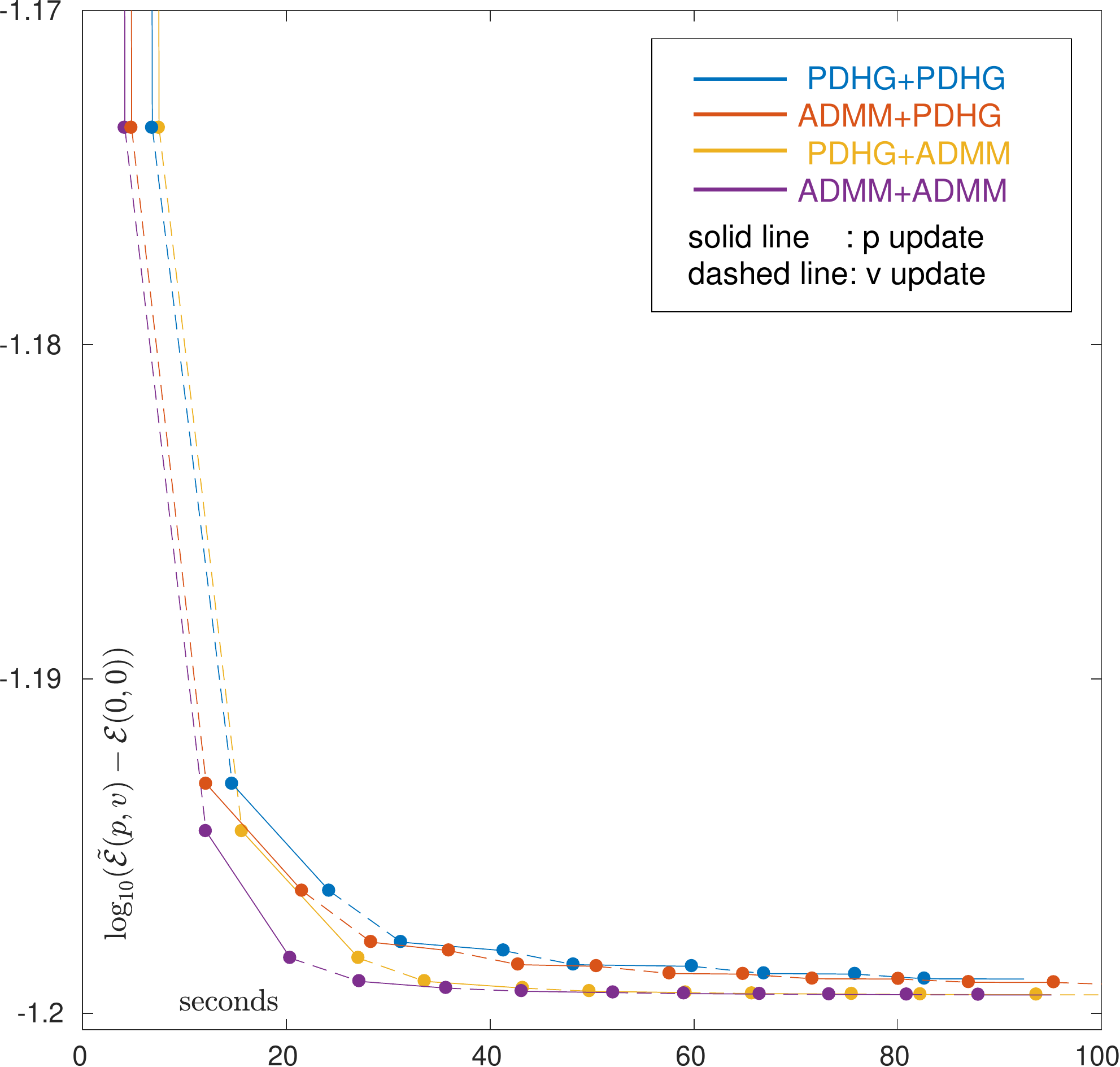}}
\hfill
\subfloat[][3D scenario\label{subfig:Cmp3D}]{\includegraphics[height=0.46\textwidth]{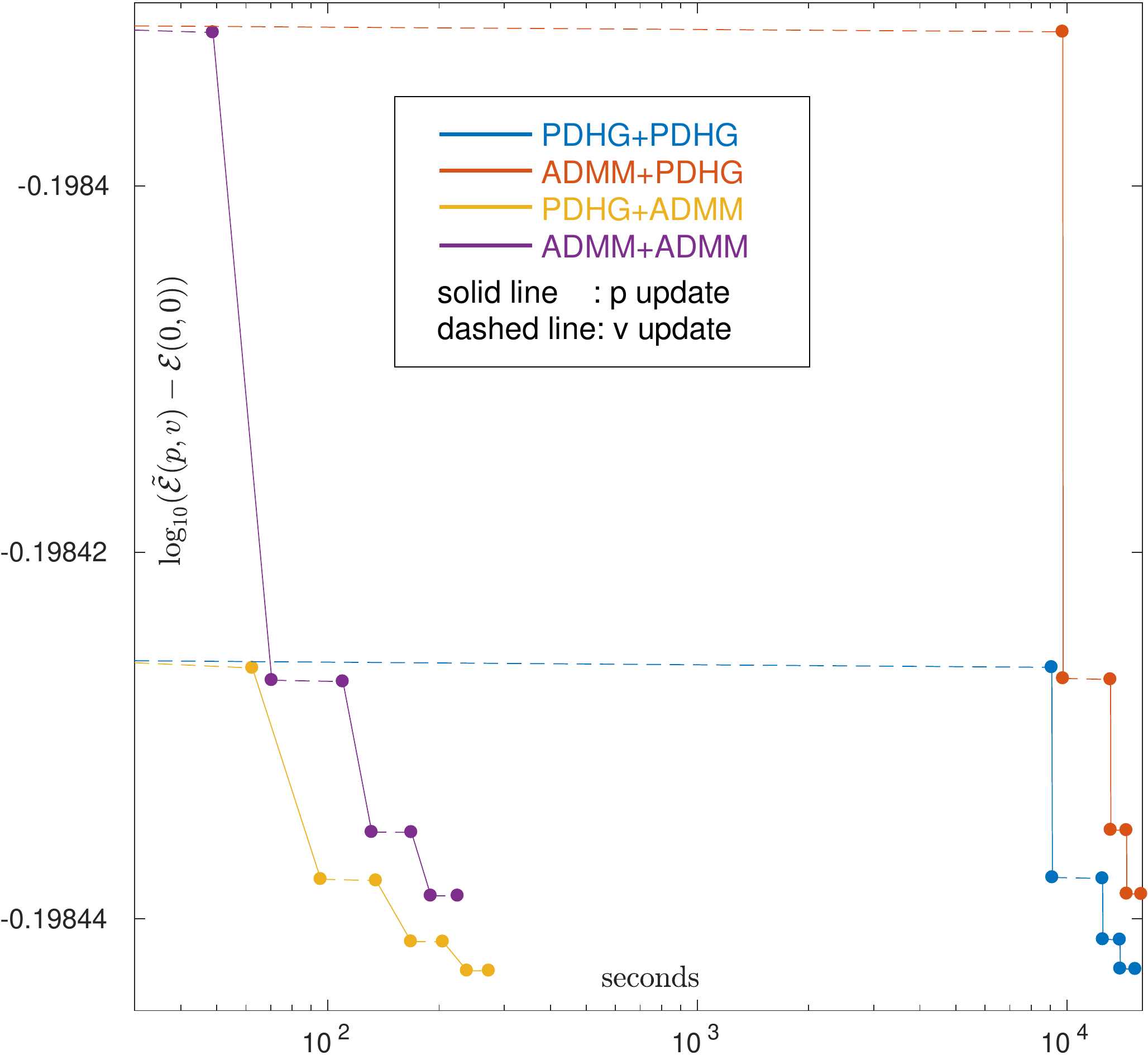}}
\caption{Comparison of different methods to solve the bi-convex optimization problem \eqref{eq:TVTVL2Denoise} via ACS for \protect\subref{subfig:Cmp2D} the 2D example described in Section \ref{sec:NumPhan} and \protect\subref{subfig:Cmp3D} $T = 10$ frames of the 3D scenario described in Section \ref{sec:ExpPhan}. The plots display the decay in energy $\tilde{\E}(p,v)$ relative to the initialization with $p = 0$, $v = 0$ vs computational time in seconds (in logarithmic scale in \protect\subref{subfig:Cmp3D}). Solid parts of the line plots correspond to the $p$ update \eqref{eq:ACSpOpt} and dashed parts to the $v$ update \eqref{eq:ACSvOpt}. A total of $4$ ACS alternations is displayed. "ADMM-PDHG" refers to using ADMM for $p$ update and PDHG for the $v$ update.}
   \label{fig:CmpOpt}
\end{figure}

\begin{figure}[tb]
   \centering
\subfloat[][\label{subfig:LS1}]{\includegraphics[height=0.48\textwidth]{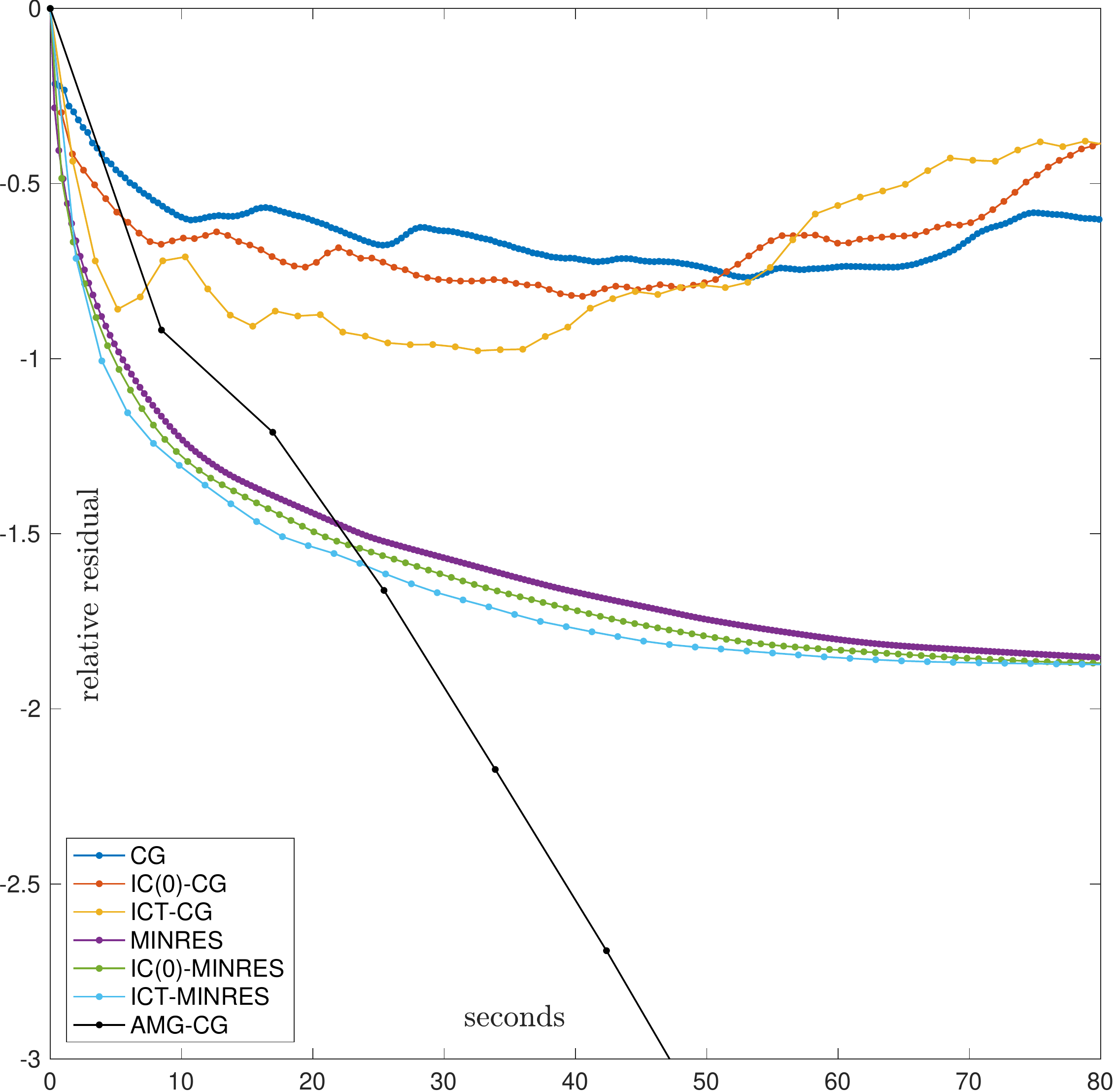}}
\hfill
\subfloat[][\label{subfig:LS2}]{\includegraphics[height=0.48\textwidth]{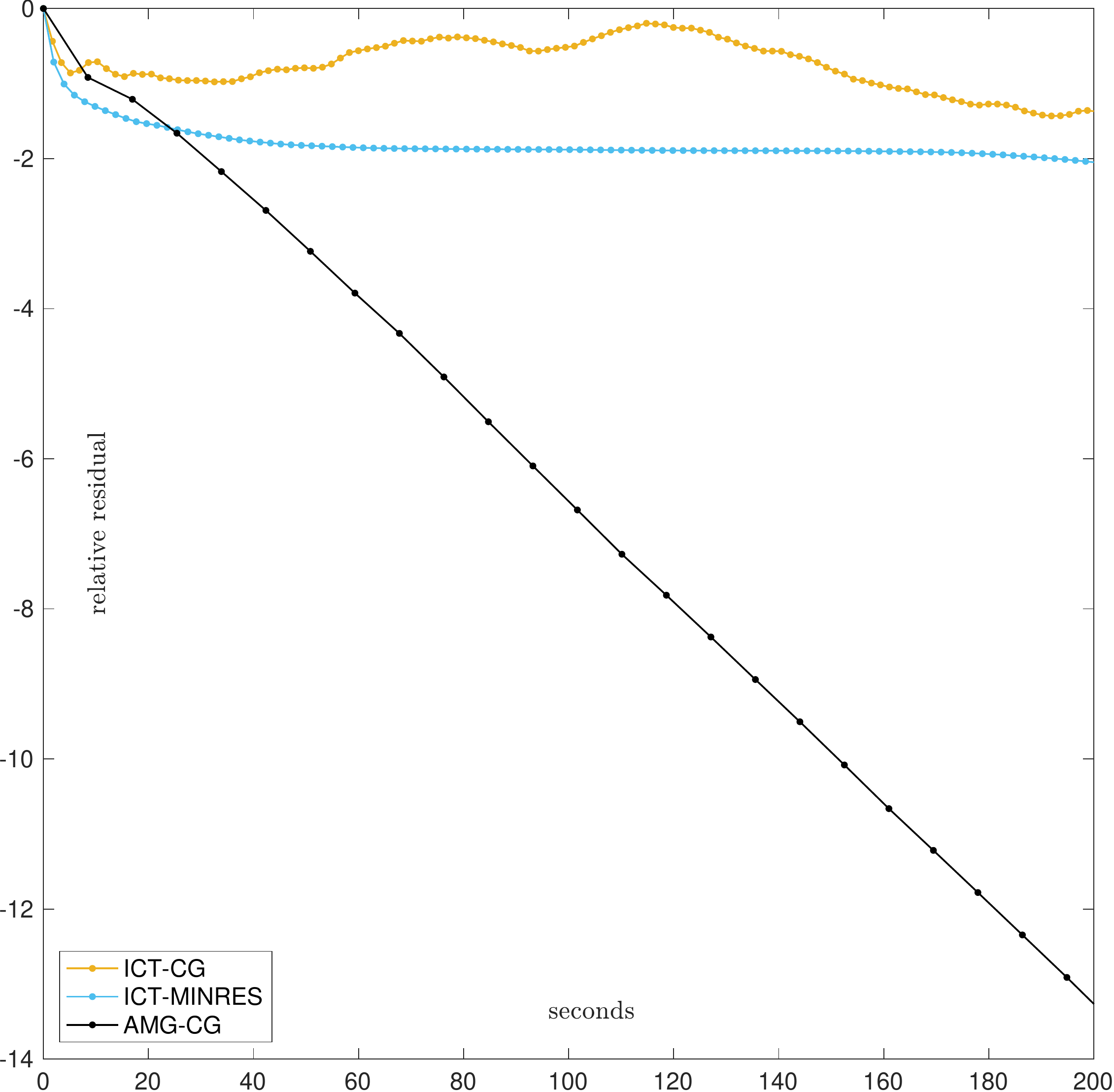}}
\caption{\protect\subref{subfig:LS1}  Comparison between different iterative methods and pre-conditioners to solve \eqref{eq:LinSys_v} (see Section \ref{sec:ExpOpt} for the details of the set-up). The vertical axis shows the relative residual while the horizontal axis shows the computation time on a single CPU core averaged over $10$ repetitions. \protect\subref{subfig:LS2} The same plot with only a subset of the solvers and expanded axis.}
   \label{fig:CmpLS}
\end{figure}

As noted earlier, the results we showed up to now were computed using ADMM in both the $p$ update \eqref{eq:ACSpOpt} and the $v$ update \eqref{eq:ACSvOpt} in the TVTVL2-regularized denoising problem \eqref{eq:TVTVL2Denoise}. In the $v$ update, AMG-CG was used as a least squares solver. In this section, we justify this choice retrospectively. Due to the large number of different parameters ACS, PDHG, and ADMM have, this is not an exhaustive comparison. We tuned all parameters we do not explicitly mention to best performance and made sure that all methods make best use of the computational platform we used (Intel Xeon CPU with 12 cores at 2.70 GHz, 256GB RAM). Another problem is caused by the non-convexity of \eqref{eq:TVTVL2Denoise} which adds an arbitrary element to such a comparison: In principle, one would need to test all methods on a large number of inputs and initilizations and compare average performances. Again, we restrict ourselves here to the two concrete examples we presented in the previous two sections and in each of those, we only examine the computation of the TVTVL2-regularized denoising problem \eqref{eq:TVTVL2Denoise} arising from the first iteration, $i = 1$, of the forward-backward splitting \eqref{eq:ProxGradA}-\eqref{eq:ProxGradB}. As $p^0 = 0$, this means we examine $\tilde{p}_t =  \nu A^T C^T f_t^c$ as an input in \eqref{eq:TVTVL2Denoise}. All other variables are initialized to $0$. Figure \ref{fig:CmpOpt} compares the decay of the denoising energy $\tilde{\E}(p,v)$ over computation time for the four different combinations of using PDHG and ADMM for each of the sub-steps. In 2D (Figure \ref{subfig:Cmp2D}), the convergence speed of the different combinations is quite similar and the different energy levels they reach corresponds to the different local minima they end up in. In 3D, the situation is quite different: Figure \ref{subfig:Cmp3D} shows that using PDHG for the $v$ update \eqref{eq:ACSvOpt} leads to prohibitively long computations times. While PDHG performs well for the $p$ update \eqref{eq:ACSpOpt} in this study, we also encountered scenarios where this is not the case. This observation was the main reason we considered using the more complicated ADMM methods in the first place: We started off by using PDHG for both sub-problems like in \cite{BuDiSc16} based on the corresponding code available on github\footnote{\href{https://github.com/HendrikMuenster/JointMotionEstimationAndImageReconstruction}{https://github.com/HendrikMuenster/JointMotionEstimationAndImageReconstruction}}. While this worked for 2D scenarios, we encountered severe difficulties for 3D scenarios which we were only able to overcome by implementing the tailored ADMM implementations presented here. \\
The main difficulty in both ADMM methods is to solve the least squares problems \eqref{eq:LinSys_p} and \eqref{eq:LinSys_v} by a fast iterative method. As explained in Section \ref{sec:ADMM}, the $v$ update \eqref{eq:LinSys_v} can be solved frame-by-frame, which allows one to explicitly set up the system matrix and use efficient pre-conditioning techniques. To compare them, we set $p^{j+1}_{t+1}$ and $p^{j+1}_{t}$ in \eqref{eq:LinSys_v} to the TV-fbf solutions shown in Figure \ref{fig:LimFBF} (note that $E = E(p_t^{j+1})$). Figure \ref{fig:CmpLS} shows the results which demonstrate that for linear systems arising from regularized 3D optical flow estimation, AMG-CG is a powerful solver.  Note, however, that this comes with increased memory costs: the system matrix is 1.24 GB large and the corresponding AMG pre-conditioner we chose here is 6.75 GB large (there is also a little computational overhead in computing them, but as they do not change over the whole ADMM scheme, this is typically negligible).

\pagebreak
\section{Discussion, Outlook and Conclusion}\label{sec:DisOutCon}
\subsection{Discussion and Outlook} \label{subsec:DisOut}
The results for both simulated and experimental data clearly demonstrate that a significant improvement of image quality over fbf reconstructions \eqref{eq:ImRecFBF} that only use spatial sparsity constraints can be obtained when using a generic spatio-temporal approach based on simultaneous, sparsity-constrained motion estimation \eqref{eq:ImRecMoEst}. Furthermore, the reconstructed motion fields provide additional information on the dynamics that can be useful for subsequent analysis. While these dynamic parameters look qualitatively correct, even from sub-sampled data, further investigations have to examine whether they are also quantitatively correct. For this first proof-of-concept study, we used very generic regularization functionals in space (TV) and a generic motion model based on a simple continuity equation \eqref{eq:OptFlowEq}. As we already obtained promising results with this rather unspecific model, we want to investigate the use of tailored motion models that better reflect the real physics of the underlying motion for a concrete application. In addition, we chose to measure the misfit to the discretized motion PDE in the squared $L_2$-norm, cf. \eqref{eq:OptFlowReg}. While this is computationally advantageous, studies generalizing this to $L_p$-norms, e.g., for $p = 1$, have shown promising results and directions for future research \cite{Di15,BuDiSc16,BuDiFrHaHeSi17}. 
\\
The main drawback of the concrete TVTVL2 model we used here \eqref{eq:TVTVL2SemiCont} is that it leads to a challenging, large-scale bi-convex optimization problem. Even with the tailored ADMM schemes we developed (cf. Section \ref{sec:ADMM}), computing the 3D reconstructions presented in Section \ref{sec:ExpPhan} took 4 days and 6 hours on a powerful work station (Intel Xeon CPU with 12 cores at 2.70 GHz, 256GB RAM, Tesla K40 GPU) compared to 6h 34m for the TV-fbf reconstruction. There are several possibilities to close this gap: 
\begin{itemize}
\item As a simple block alternation, the ACS scheme can be modified by introducing techniques like over-relaxation, inertia methods or line-search. 
\item For solving sub-step \eqref{eq:ACSvOpt}, developing an ADMM scheme that uses an algebraic multigrid pre-conditioner was crucial. However, the high memory demand of this approach limits the number of frames which can be computed in parallel. Using \emph{geometric multigrid} pre-conditioning instead could keep the fast convergence (cf. Figure \ref{fig:CmpLS}) while requiring much less memory \cite{BrWeKoSc06}.
\item If the non-smooth sparsity constraints are approximated by smooth functionals such as the Huber functional, fast, monotone solvers can be used to solve \eqref{eq:ACSpOpt} and \eqref{eq:ACSvOpt}, see, e.g. \cite{Vo02}.
\end{itemize}
Another potential problem is that the ACS scheme presented in Section \ref{sec:BiConv} will only converge to a local minimum of the bi-convex variational energy \eqref{eq:TVTVL2SemiCont}. Figure \ref{fig:CmpOpt} showed that already the choice of the convex optimization scheme  to solve \eqref{eq:ACSpOpt} and \eqref{eq:ACSvOpt} can influence which local minimum is found. Other parameters like the accuracy with which these problems are solved, how the schemes are initialized, whether the scenario is 2D or 3D, etc., have an often non-trivial influence as well. In future work, we plan to examine these issues in a systematic way.\\
From a modelling perspective, the simple optical flow discretization \eqref{eq:OptFlowReg} we chose here can only resolve small motions: If the support of $p_{t+1} - p_t$ and $\NabCD p_t$ do not overlap, $v_t$ cannot minimize $\sqnorm{p_{t+1} - p_t + (\nabla p_t) \cdot v_t}$. An extension of the framework to estimate large-scale motions is described in \cite{Di16}.\\
This article focused on the mathematical and computational aspects of 4D PAT. We therefore assumed here that there is a generic binning of the sequence of acoustic measurements \eqref{eq:GenScan} into temporal bins during which the target can be considered static (and used phantoms for which this holds true) and only compared image quality for a fixed sub-sampling factor. However, in reality, sequential scanners measure a single time pressure course for every pulse of the excitation laser. The temporal binning of this stream of acquisitions leads to a more complicated interplay between artefacts arising from sub-sampling, motion-blur and the spatio-temporal continuity imposed by the variational model. We will examine this issue more closely in forthcoming work that will focus on the technical and practical aspects of 4D PAT with novel acoustic scanners \cite{HuOgZhCoBe16,NaLuZhBeArBeCo17}. For the application to \emph{in-vivo} imaging, additional challenges need to be addressed, such as heterogeneous tissue properties.

\subsection{Conclusion}\label{subsec:Con}

In this work, we extended our earlier results on using compressed sensing techniques to accelerate high resolution 3D PAT acquisition with sequential scanners \cite{ArBeBeCoHuLuOgZh16}. We demonstrated that in the context of dynamic PAT, another substantial increase of image quality can be obtained by using a generic variational framework that couples sparsity-constrained image reconstruction and simultaneous, sparsity-constrained motion estimation. In particular, we considered a motion model based on the popular optical flow equation and used the total variation functional as sparsity constraints. For this, promising results for simulated and experimental data were obtained in a proof-of-concept study that justifies further research in this field. A major challenge for using these variational approaches for large scale 4D inverse problems with complicated forward operators are the computational demands of the corresponding optimization routines. We described and examined a set of related methods that can be used as a starting point to implement similar strategies for other applications.  


\clearpage

\appendix

\section{Proximal Operators} \label{sec:ProxOp}

An extensive overview on how to use and compute proximal operators \eqref{eq:Prox} is given in \cite{CoPe11}. The splits we use in this work have been introduced such that the functionals for which we have to compute the proximal operators decouple over space and time into the sum of 1 or $d$ dimensional functionals $\phi(x)$ or $\phi(x_1, \ldots, x_d)$. As such, all proximal operators can be computed explicitly and point-wise in space and time, i.e., for an image/vector field sequence $x \in \R^{N T}$/$x \in \R^{d N T}$, the proximal operators can be computed by solving $N T$ sub-problems of dimension 1/$d$ using explicit formulae. \\
For $\mathcal{G}(x)$ in \eqref{eq:G_PDHG_p}  this leads to 
\begin{equation}
\phi(x) = \chi_+(x) + (x - z)^2\enspace, \qquad 	\prox_{\alpha \phi}(\tilde{x}) = \max \left( 0,  \alpha z  + \frac{\tilde{x}}{\alpha +1}\right) \enspace .
\end{equation}
The proximal operator for the functional $\mathcal{G}(x)$ in \eqref{eq:G_PDHG_v} is a $d$-dimensional quadratic problem:
\begin{align}
\phi(x) &=  \frac{1}{2} (z + c_1 x_1 + \ldots + c_d x_d)^2 \\ 
\prox_{\alpha \phi}(\tilde{x}) &= \argmaxsub{x \in \R^d} \left\lbrace \frac{\alpha}{2} \left(z + \sum_i^d c_i x_i\right)^2  + \frac{1}{2} \sum_i^d \left(x_i  -\tilde{x}_i\right)^2 \right\rbrace
\end{align}
Its optimality condition leads to a $d$-dim linear system, which we show here for $d=3$:
\begin{equation}
\quad \begin{bmatrix}
(1+ \alpha c_1^2) & \alpha c_1 c_2    & \alpha c_1 c_3   \\
\alpha c_2 c_1    & (1+ \alpha c_2^2) & \alpha c_2 c_3   \\
\alpha c_3 c_1    &  \alpha c_3 c_2   & (1+ \alpha c_3^2)\\
\end{bmatrix} 
\begin{bmatrix}
x_1 \\
x_2 \\
x_3 \\
\end{bmatrix} 
= 
\begin{bmatrix}
\tilde{x}_1 - \alpha c_1 z \\
\tilde{x}_2 - \alpha c_2 z \\
\tilde{x}_3 - \alpha c_3 z \\
\end{bmatrix} 
\end{equation}
It can be solved explicitly for $d= 2, 3$ and for its use within the PDHG scheme, most relevant terms can be precomputed. \\
The $\ell_1$-norms involved in the isotropic TV terms are actually global $\ell_1$ norms of the local $\ell_{2}$ norms of the gradient vectors. For a gradient field of image $z$ represented as $y \in R^{d N}$ indexed as $y_{{x_i},j}$ for derivative direction and location index, respectively, we have 
\begin{align}
\norm{\NabFD z}_1 &= \sum_j^N \sqrt{\sum_i^d y_{{x_i},j}^2} \enspace, \qquad \Longrightarrow \qquad \phi(x) = \sqrt{x_1 + \ldots + x_d} \enspace, \\
\prox_{\alpha \phi}(\tilde{x}) &= 
\begin{cases}
\max(\phi(\tilde{x}) - \alpha , 0)  \; \tilde{x}/\phi(\tilde{x}) \hspace{1em} \text{if} \hspace{1em} \phi(\tilde{x}) > 0 \\
0 \hspace{11.2em} \text{else}
\end{cases}
\enspace .
\end{align}
With this, one can easily build the proximal operator for \eqref{eq:F_PDHG_v} and \eqref{eq:F_ADMM_p}. \\
Next we need the convex conjugates $\mathcal{F}^*(y)$ and their proximal mappings in some places. For the isotropic TV term, we have 
\begin{equation}
\mathcal{F}(y) = \alpha \norm{y}_1 = \alpha \sum_j^N \sqrt{\sum_i^d y_{{x_i},j}^2} \;,  \quad \mathcal{F}^*(y) =  \sum_j^N \alpha \chi_{[0,1]}\left( \frac{1}{\alpha} \sqrt{\sum_i^d y_{{x_i},j}^2} \right) \;,
\end{equation}
which means that $\mathcal{F}^*(y)$ is $0$ if all gradient vectors have an amplitude that is smaller than $\alpha$ and $\infty$ else, see, e.g., \cite{ChPo11}. As such, the proximal operator is just a projection:
\begin{equation}
\phi(y) = \chi_{[0,1]} \left( \frac{1}{\alpha} \sqrt{\sum_i^d y_{x_i}^2} \right) \quad \Longrightarrow \quad \prox_{\beta \phi}(\tilde{y}) = \frac{\tilde{y}}{\max \left(1, \frac{1}{\alpha} \sqrt{\sum_i^d \tilde{y}_{x_i}^2} \right)}
\end{equation}
The second part of $\mathcal{F}(y)$ in \eqref{eq:G_PDHG_v} is $\frac{\tilde{\gamma}}{2}\sqnorm{y_2}$. Its convex conjugate is given by $\frac{1}{2 \tilde{\gamma}}\sqnorm{y_2}$ and the proximal mapping can be computed using
\begin{equation}
\phi(y) = \frac{y^2}{2 \tilde{\gamma}} \quad \Longrightarrow \quad \prox_{\alpha \phi}(\tilde{y}) = \frac{\tilde{\gamma}}{\tilde{\gamma} + \alpha} \; \tilde{y} \enspace .
\end{equation}


\section*{Acknowledgments}
We would like to thank Hendrik Dirks for very helpful discussions and support for his code\footnote{\url{https://github.com/HendrikMuenster/JointMotionEstimationAndImageReconstruction}} which provided a template for our implementation of the TVTVL2-denoising function. Further more, we gratefully acknowledge the support of NVIDIA Corporation with the donation of the Tesla K40 GPU used for this research.

\bibliographystyle{siamplain}
\bibliography{all}
\end{document}